\documentclass[12pt]{article}

\scrollmode
\usepackage[latin1]{inputenc}
\usepackage{amsmath,amssymb,amsfonts}

\newtheorem{thm}{Theorem}[section]
\newtheorem{prop}[thm]{Proposition}
\newtheorem{lem}[thm]{Lemma}

\newtheorem{defn}[thm]{Definition} 
\newtheorem{exmp}[thm]{Example}   
\newtheorem{rem}[thm]{Remark}   
\newtheorem{cor}[thm]{Corollary}

\def\be#1 {\begin{equation} \label{#1}}
\newcommand{\ee}{\end{equation}}

\def\sqw{\hbox{\rlap{\leavevmode\raise.3ex\hbox{$\sqcap$}}$%
\sqcup$}}
\def\findem{\ifmmode\sqw\else{\ifhmode\unskip\fi\nobreak\hfil
\penalty50\hskip1em\null\nobreak\hfil\sqw
\parfillskip=0pt\finalhyphendemerits=0\endgraf}\fi}

\newcommand{\mb}{\medskip\noindent}
\newcommand{\gb}{\bigskip\noindent}
\newcommand{\R}{\mathbb R}

\newcommand{\N}{\mathbb N}

\newcommand{\C}{\mathbb C}
\newcommand{\s}{\mathbf S}

\def\deme {\noindent {\bf Proof : }}

\title{New Abstract Hardy Spaces}
\author{Fr\'ed\'eric Bernicot and Jiman Zhao
\\
 \small Universit\'e de Paris-Sud, Orsay et CNRS 8628, 91405 Orsay Cedex, France
\\
\small {\em E-mail address:} {Frederic.Bernicot@math.u-psud.fr}\\
\small School of Mathematical Sciences,
Beijing Normal University,\\
\small Key Laboratory of Mathematics and Complex Systems, Ministry
of Education,\\
\small Beijing 100875,
P.R. China\\
\small {\em E-mail address:} {jzhao@bnu.edu.cn}\\
}

\begin{document}

\maketitle

\begin{abstract}
  The aim of this paper is to propose an abstract construction of spaces which
  keep the main properties of the (already known) Hardy spaces $H^1$. We construct
  spaces through an atomic (or molecular) decomposition. We prove some
  results about continuity from these spaces into $L^1$ and some results about
  interpolation between these spaces and the Lebesgue spaces.
  We also obtain some results on weighted norm inequalities.
  Then we apply this abstract theory to the $L^p$ maximal regularity. Finally we present partial results in order to understand a characterization of the duals of Hardy spaces.

\vspace{0.1in} \footnotesize {\bf Key words}: Hardy spaces, atomic decomposition, interpolation, maximal regularity \\
\vspace{0.1in}
 {\bf  AMS2000 Classification}: 42B25, 42B30, 42B35, 46M35

\end{abstract}

\section{Introduction}

The theory of real Hardy spaces started in the 60's, and in the 70's
the atomic Hardy space appeared. Let us recall its definition first
(see \cite{CW}).

\gb Let $(X,d,\mu)$ be a space of homogeneous type. Let $\epsilon>0$
be a fixed parameter. A function $m\in L^{1}_{loc}(X)$ is called an
$\epsilon$-molecule associated to a ball $Q$ if $\int_X m d\mu=0$,
for all $i\geq 0,$
$$ \left( \int_{2^{i+1}Q \setminus 2^i Q} |m|^2 d\mu \right)^{1/2} \leq \mu(2^{i+1}Q)^{-1/2} 2^{-\epsilon i} \textrm{  and  } \left( \int_Q |m|^2 d\mu \right)^{1/2} \leq \mu(Q)^{-1/2}$$
We call $m$ an atom if in addition we have $ \textrm{supp} (m) \subset Q$. So an atom is exactly an $\infty$-molecule.
Then a function $f$ belongs to $H^1_{CW}(X)$ if there exists a decomposition
$$f =\sum_{i\in\N} \lambda_i m_i  \qquad \mu-a.e, $$
where $m_i$ are $\epsilon$-molecules and $\lambda_i$ are
coefficients which satisfy
$$\sum_{i} |\lambda_i| <\infty.$$
It was proved in \cite{steine} that the whole space $H^1_{CW}$ does not depend on $\epsilon$, as in fact one obtains the same space replacing $\epsilon$-molecules by atoms or $\epsilon'$-molecules with $\epsilon'>0$.

\gb In the Euclidean case ($X=\R^n$ with the Lebesgue measure) this
space has many different characterizations, thanks to \cite{steine}~:
\begin{align}
 f\in H^1_{CW}(\R^n)  & \Longleftrightarrow f\in {\mathcal H}^{1} := \left\{f\in L^1(\R^n)\, ;\, \nabla (\sqrt{-\Delta})^{-1}(f) \in L^1(\R^n,\R^n) \right\} \label{Hrond} \\
   & \Longleftrightarrow x \mapsto \sup_{\genfrac{}{}{0pt}{}{y\in\R^n,t>0}{|x-y|\leq t}} \left| e^{-t\sqrt{-\Delta}}(f)(y) \right| \in L^1(\R^n) \label{cha1} \\
  & \Longleftrightarrow x \mapsto \left(\int_{\genfrac{}{}{0pt}{}{y\in\R^n,t>0}{|x-y|\leq t}} \left| t\nabla e^{-t\sqrt{-\Delta}}(f)(y) \right|^2 \frac{dydt}{t^{n+1}} \right)^{1/2} \in L^1(\R^n), \label{cha2}
\end{align}
where $\nabla (\sqrt{-\Delta})^{-1}$ is the Riesz transform. The
space ${\mathcal H}^1$ defined by (\ref{Hrond}) was the original
Hardy space of E.M. Stein (see \cite{Se}) and \cite{steine} provided
the equivalence with the definition using the maximal function and the area
integral. The link with $H^1_{CW}(\R^n)$ (due to R. Coifman one year later
in \cite{FS}) can be understood from the celebrated theorem of C.
Fefferman which says (in vague terms)
$$ h \in ({\mathcal H}^1)^{*} \Longleftrightarrow \|h\|_{BMO}:= \left(\sup_{Q \textrm{ ball}} \frac{1}{\mu(Q)}\int_Q \left|f - \frac{1}{\mu(Q)}\int_Q f d\mu\right|^2 d\mu\right)^{1/2} <\infty.$$
Here $BMO$ is the space of John-Nirenberg. In fact it is relatively easy to show that
$$ h \in (H_{CW}^1(\R^n))^{*} \Longleftrightarrow \|h\|_{BMO} <\infty,$$
hence the identification between ${\mathcal H}^1$ and $H^1_{CW}(\R^n)$.

\gb The space $H^1_{CW}(X)$ is a good substitute of $L^1(X)$ for many reasons. For instance, Calder\'on-Zygmund operators map $H^1_{CW}(X)$ to $L^1(X)$ whereas they do not map $L^1(X)$ to $L^1(X)$. In addition, $H^1_{CW}(X)$ (and its dual) interpolates with Lebesgue spaces $L^p(X)$, $1<p<\infty$. That is why $H^1_{CW}(X)$ is a good space to extend the scale of Lebesgue spaces $(L^p(X))_{1<p<\infty}$ for $p$ tends to $1$ and its dual $BMO$ for $p$ tends to $\infty$.

\gb However, there are situations where $H^1_{CW}(X)$ is not the
 right substitute to $L^1(X)$ and there has been recently a number of works
 with goal to define an adapted Hardy space \cite{A1,A,AR,D,DY,D1,D2,DY1,DZ,HM}.
 For example $H^1_{CW}(X)$ is not well adapted to operators such as the Riesz transform
  on Riemannian manifolds (or on graphs) or the maximal regularity
  operator \cite{A1,A,ACDH,HM,R}. That is why in \cite{DY,HM} the authors
   define a new space $H^1_{L}$ by the $L^1(X)$ norm of the previous maximal functions
   (in (\ref{cha1}) and (\ref{cha2})) with another operator $L$ instead of the
    Laplacian $\Delta$. With this new definition, they show that under some conditions on
    $L$, the space $H^1_{L}$ has an equivalent molecular definition. They have studied the intermediate spaces
  between $H^1_{L}$ and the Lebesgue spaces and also the dual space $(H^1_{L})^{*}$.

\gb Our aim here is to construct abstract Hardy spaces by a molecular (or atomic) decomposition
 and we want to use the weakest assumptions to obtain good properties for these spaces.
 Mainly, we want to have a criterion for the continuity of an operator from the Hardy
 space into $L^1(X)$ and to have an interpolation result between the Hardy space and Lebesgue
 spaces.
 Then we will apply these abstract results to a particular case and we will obtain
 results about maximal $L^q$ regularity. We will finish this paper by
  the study of the dual space.

\section{Definitions}

\gb Let $(X,d,\mu)$ be a space of homogeneous type. Excepted for the section \ref{seven} where the space is important, we omit the space $X$ and we shall write $L^p$ for $L^p(X,\R)$ if no confusion arises. Here we are working with real valued function and we will use "real" duality. We have the same results with complex duality and complex valued functions. \\
 So $d$ is a
quasi-distance on the space $X$ and $\mu$ a Borel measure which
satisfies the doubling property~: \be{homogene} \exists A>0,  \
\exists \delta>0, \qquad \forall x\in X,\forall r>0, \forall t\geq
1,\qquad \frac{\mu(B(x,tr))}{\mu(B(x,r))} \leq A t^{\delta}, \ee
where $B(x,r)$ is the open ball with center $x\in X$ and radius
$r>0$. We call $\delta$ the homogeneous dimension of $X$. For $Q$ a ball,
and $i\geq 0$, we write $S_i(Q)$ the scaled
corona around the ball $Q$~:
$$ S_i(Q):=\left\{ x,\ 2^{i} \leq 1+\frac{d(x,c(Q))}{r_Q} < 2^{i+1} \right\},$$
where $r_Q$  is the radius of the ball $Q$ and $c(Q)$ its center.
 Note that $S_0(Q)$ corresponds to the ball $Q$ and $S_i(Q) \subset 2^{i+1}Q$ for $i\geq 1$, where $\lambda Q$ is as usual the ball with center $c(Q)$ and radius $\lambda r_Q$.

\mb Let us denote by $\mathcal{Q}$ the collection of all balls~:
$$ \mathcal{Q}:= \left\{ B(x,r),\ x\in X, r>0 \right\}.$$
Let $\mathbb{B}:=(B_Q)_{Q\in \mathcal{Q}}$ be a collection of $L^2$-bounded linear operator, indexed by the collection $\mathcal{Q}$.
We assume that these operators $B_Q$ are uniformly bounded on $L^2$ : there exists a
constant $0<A'<\infty$
so that~: \be{operh} \forall f\in L^2 ,\ \forall  Q \textrm{ ball}, \qquad \|B_Q(f) \|_2 \leq A'\|f\|_2. \ee

\mb In the rest of the paper, we allow the constants to depend on $A$, $A'$ and $\delta$.

\mb Now, we define  {\it atoms} and {\it molecules} by using the
collection $\mathbb{B}$~:

\begin{defn} Let $\epsilon>0$ be a fixed parameter. A function $m\in L^{1}_{loc}$ is called an $\epsilon$-molecule associated to a ball $Q$ if there exists a real function $f_Q$ such that
$$m=B_Q(f_Q),$$
with
$$\forall i\geq 0, \qquad  \|f_Q\|_{2,S_i(Q)} \leq \left(\mu(2^{i}Q)\right)^{-1/2} 2^{-\epsilon i}.$$
We call $m=B_Q(f_Q)$ an atom if in addition we have $supp(f_Q) \subset Q$. So an atom is exactly an $\infty$-molecule.
\end{defn}

\mb The functions $f_Q$ in this definition are normalized in
$L^1$. It is easy to show that
$$ \|f_Q\|_{1} \lesssim 1 \qquad \textrm{and} \qquad \|f_Q\|_{2} \lesssim \mu(Q)^{-1/2}.$$
So by the $L^2$-boundedness of the operator $B_Q$,
we have that each molecule belongs to the space $L^2$. However a molecule is
not (for the moment) in the space $L^1$. In Section \ref{particular},
 we will put some further conditions on the operators $B_Q$ which will guarantee that
  the molecules will form a bounded set in the space $L^1$. But for the moment,
  we want to work with the most general case. Now we are able to define our abstract
  Hardy spaces~:

\begin{defn} A measurable function $h$ belongs to the molecular
Hardy space $H^1_{\epsilon,mol}$ if there exists a decomposition~:
$$h=\sum_{i\in\N} \lambda_i m_i  \qquad \mu-a.e, $$
where for all $i$, $m_i$ is an $\epsilon$-molecule and $\lambda_{i}$
are real numbers satisfying
$$\sum_{i\in \N} |\lambda_i| <\infty. $$
We define the norm~:
$$\|h\|_{H^1_{\epsilon,mol}}:= \inf_{h=\sum_{i\in\N} \lambda_i m_i} \sum_{i} |\lambda_i|.$$
Similarly we define the atomic space $H^1_{ato}$ replacing $\epsilon$-molecules by atoms.
\end{defn}

\mb Let us make some remarks.

\begin{rem}

\mb $1-)$ First we only ask that the decomposition
$$h(x) =\sum_{i\in\N} \lambda_i m_i(x) $$
is well defined for almost every $x\in X$. So the assumption is very
weak and it is possible that the measurable function $h$ does not
belong to $L^1_{loc}$. It is not clear whether these abstract normed vector
spaces are complete. The problem is that we do not know whether
the
decompositions for $h$ converge absolutely. See Section \ref{particular}
 for a condition insuring this. However  we do not need completeness for the moment. \\
$2-)$ We have the following continuous inclusions~:
\be{inclusio} \forall\, 0<\epsilon<\epsilon',\  \qquad  H^1_{ato} \hookrightarrow H^1_{\epsilon',mol} \hookrightarrow H^1_{\epsilon,mol}. \ee
In fact the space $H^1_{ato}$ corresponds to the space $H^1_{\infty,mol}$.
 For $0<\epsilon<\epsilon'<\infty$ the space $H^1_{\epsilon',mol}$ is
 dense in $H^1_{\epsilon, mol}$. In the general case,
 it seems to be very difficult to study the dependence of the space $H^1_{\epsilon, mol}$ with the parameter $\epsilon$ and we will not study this question here.
 \\
$3-)$ We have seen that each molecule is an $L^2$ function. So it is obvious that $L^2 \cap H^1_{\epsilon, mol}$ is dense in $H^1_{\epsilon, mol}$ and that $L^2 \cap H^1_{ato}$ is dense in $H^1_{ato}$.
\end{rem}

\section{Comparison with other Hardy spaces.} \label{sectionexample}

\subsection{The space of Coifman-Weiss.}

\gb Due to its atomic definition, the Hardy space $H^1_{CW}(\R^n)$ of Coifman-Weiss is obtained by choosing the operator $B_Q$ as follows~:
$$ B_Q(f)(x)=f(x){\bf 1}_{Q}(x) - |Q|^{-1}\left(\int_Q f\right) {\bf 1}_{Q}(x).$$
Our atoms are the same as the ones defined in \cite{steine}.
However, our  molecule is different from the one in \cite{steine}.
In fact, because the $B_Q$ has the property that supp $
B_Q(f) \subset Q$ for any $f$, atoms and $\epsilon$-molecules are
the same with our definition. Hence, for $B_Q$ with this specific
property, our atomic and molecular spaces are the same.

\subsection{Hardy spaces for Schrödinger operators with potentials.}

\mb Let $X=\R^{n}$ and $V$ a nonnegative function on $X$. We consider the Schrödinger operator
$$ L(f)(x):=-\Delta f(x) + V(x) f(x).$$

\mb {\bf First case :} $V$ is compactly supported and $V\in L^p$ with $2p>n\geq 3$ (we refer to \cite{DZ} for the details). \\
By this assumption, it is well known that $-L$ generates a
$L^2$-bounded semigroup $(K_t)_{t>0}$, whose kernels satisfy some
gaussian estimates. J. Dziubanski and J. Zienkiewicz
define a Hardy space $H^1_{L}$ by a maximal operator. A
function $f\in L^1(\R^n)$ belongs to $H^1_{L}$ if
$$ \|f\|_{H^1_L} := \left\| \sup_{t>0} \left|K_t f(x) \right| \right\|_{1} <\infty.$$
Using the properties of semi-group, the authors introduce
$$\omega(x):= \lim_{s\to \infty} K_s{\bf 1}_{\R^n}(x)$$
 and prove that the limit exists and that there is a constant $c>0$ such that
\be{fonome} \forall x\in \R^n, \qquad c \leq \omega(x) \leq 1.\ee
With this function, they obtain an atomic decomposition of their
Hardy space $H^1_{L}$ with the following definition of atoms~: a
function $b$ is a $H^1_{L}$-atom if there exists a ball $Q$ such
that $b$ is supported in $Q$ and satisfies
$$ \|b\|_{2} \leq |Q|^{-1/2} \qquad \textrm{and} \qquad \int_Q b(x) \omega(x)dx =0.$$
So we can identify their space $H^1_{L}$ with our atomic space $H^1_{ato}$ by choosing our operators $B_Q$ as
$$ B_Q(f)(x) := f(x){\bf 1}_{Q}(x) -\left( \frac{1}{\omega(Q)} \int f(x) \omega(x)dx \right){\bf 1}_{Q}(x).$$
Then due to (\ref{fonome}), the two definitions of atoms are equivalent and so
$$ H^1_{L} = H^1_{ato}.$$

\mb {\bf Second case :} $V$ is a nonnegative polynomial (we refer to \cite{D,D1} for
the details). \\
In this case, it is well known that $-L$ generates a $L^2$-bounded
semigroup $(T_t)_{t>0}$, which satisfies some gaussian estimates. J.
Dziubanski defines a Hardy space in the same way as above : a
function $f\in L^1$ belongs to $H^1_{L}$ if
$$ \|f\|_{H^1_L} := \left\| \sup_{t>0} \left|T_t f(x) \right| \right\|_{1} <\infty.$$
Let
$$m(x,V):= \sum_{\beta} \left| D^\beta V(x) \right|^{1/(|\beta|+2)}$$
which is bounded below by a constant $c>0$. In \cite{D1} the author shows an atomic decomposition of this space with the following definition~: a function $b$ is an $H^1_{L}$-atom if there exists a ball $Q=B(y_0,r)$ with
$$ \textrm{supp}(b) \subset Q, \quad \|b\|_{2} \leq |Q|^{-1/2},\quad r\leq m(y_0,V)^{-1} $$
and if $r\leq \frac{1}{4} m(y_0,V)^{-1}$ then
$$ \qquad \int_{Q} b(x) dx =0 . $$
This definition of atoms is a particular case of ours if we define the operator $B_Q$ for $Q=B(y_0,r)$ a ball by~:
$$ B_Q(f)(x) := \left\{ \begin{array}{l}
 f(x){\bf 1}_{Q}(x)  \qquad \textrm{  if $\frac{1}{4}m(y_0,V)^{-1}<r\leq m(y_0,V)^{-1}$} \\
 f(x){\bf 1}_{Q}(x)-\left(\frac{1}{|Q|} \int_Q f(y)dy \right){\bf 1}_{Q}(x) \qquad  \textrm{  if $r\leq \frac{1}{4} m(y_0,V)^{-1}$} \\
 0 \qquad \textrm{  if $r> m(y_0,V)^{-1}$.}
\end{array}
\right. $$
With this choice we have
$$ H^1_{L} = H^1_{ato}.$$
It is shown in \cite{D2} that one can take $w$ as a reverse Hölder weight with exponent $n/2$ to obtain an identical atomic decomposition with the measure $\omega(x)dx$ instead of the Lebesgue measure.

\mb

\subsection{Hardy spaces associated to divergence form elliptic operators.}
\label{elliptic}

\mb Let $X=\R^n$ and $A$ be an $n\times n$ matrix-valued function satisfying the ellipticity condition~: there exist two constants $\Lambda\geq \lambda>0$ such that
$$ \forall \xi,\zeta\in \C^n, \qquad \lambda|\xi|^2 \leq Re \left( A\xi \cdot \overline{\xi} \right) \quad  \textrm{and} \quad |A\xi \cdot \overline{\zeta} | \leq \Lambda|\xi||\zeta|.$$ We define the second order divergence form operator
$$ L(f):= -\textrm{div} (A \nabla f).$$
In \cite{HM} S. Hofmann and S. Mayboroda define a Hardy space $H^1_{L}$ associated to
this operator and give several charaterizations. For $f\in L^1$ we have
the equivalence of the following norms~:
\begin{align*}
 \|f\|_{H^1_{L}} & := \|f\|_{1}+ \left\| \left(\iint_{\genfrac{}{}{0pt}{}{t>0,\ y\in\R^n}{|x-y|\leq t}} \left|t^2L e^{-t^2L} f(y) \right|^2 \frac{dtdy}{t^{n+1}} \right)^{1/2} \right\|_{1} \label{maximaldef} \\
  & \simeq \|f\|_{1} + \left\| \sup_{\genfrac{}{}{0pt}{}{t>0,\ y\in\R^n}{|x-y|\leq t}} \left(\frac{1}{t^n} \int_{B(y,t)} \left|e^{-t^2L} f(z) \right|^2 dz \right)^{1/2} \right\|_{1}. \nonumber
  \end{align*}
In addition, they prove a molecular decomposition with the following definition : let $\epsilon>0$ and $M>n/4$ be fixed, a function $m\in L^2$ is a $H^1_{L}$-molecule if there exists a ball $Q\subset \R^n$ such that~:
\begin{align}
& \forall i\geq 0, \qquad \|m\|_{2,S_i(Q)} \leq 2^{-i\epsilon} |2^{i+1}Q|^{-1/2} & \\
& \forall i\geq 0, \forall k\in\{1,...,M\}, \qquad \left\| \left(r_Q^{-2}L^{-1}\right)^{k} m \right\|_{2,S_i(Q)} \leq 2^{-i\epsilon} |2^{i+1}Q|^{-1/2}. &
\end{align}
We do not know how to realize these molecules with our definition. However with
$$ B_Q(f) := (r_Q^2L)^{M}e^{-r_Q^2L}(f) \textrm{  or  }  B_Q(f) := \left(Id - (Id+r_Q^2L)^{-1} \right)^{M}(f),$$
our $\epsilon$-molecules are $H^1_{L}$-molecules. So we have the
inclusion
$$ H^1_{\epsilon, mol} \hookrightarrow H^1_{L}.$$

\subsection{Hardy spaces associated to a general semigroup.}
\label{elliptic2}

\gb Let $X =\R^n$. In \cite{DY}, X.T. Duong and L. Yan have study
the space $H^1_{L}$ with a more general operator $L$. They assume
that there exists $\omega\in (0,\pi/2)$ such that $L$ generates a
holomorphic semigroup $e^{-zL}$ with $0\leq |Arg(z)|<\pi/2-\omega$,
which has a $H_\infty$ calculus on $L^2(\R^n)$ and gaussian estimates for its kernel. Then they
define  a Hardy space $H^1_{L}$ by : for all functions $f\in L^1$
 $$ f\in H^1_{L} \Longleftrightarrow \left\| \left(\iint_{\genfrac{}{}{0pt}{}{t>0,\ y\in\R^n}{|x-y|\leq t}} \left|t^2L e^{-t^2L} f(y) \right|^2 \frac{dtdy}{t^{n+1}} \right)^{1/2} \right\|_{1}.$$
They obtained a molecular decomposition using tent spaces
(Proposition 4.2 of \cite{DY}) : a function $m$ is called a
$H^1_{L}$-molecule associated to a ball $Q \subset \R^n$ if
$$ m(x) = \int_0^\infty  t^2L e^{-t^2L}[a(t,.)](x) \frac{dt}{t}$$
with $a(t,x)$ satisfying
$$ \textrm{supp}(a) \subset \left\{ (x,t)\in X\times (0,\infty),\ x\in Q, 0<t\leq r_Q\right\} \quad \textrm{and} \quad \|a\|_{2,\frac{dxdt}{t}} \leq |Q|^{-1/2}.$$
This definition of molecules is probably less restrictive than ours. So with the choice
$$ B_Q(f):= (r_Q^2L)e^{-r_Q^2L}(f) \textrm{  or  } B_Q(f):= f-e^{-r_Q^2L}(f),$$
we only know the inclusion
$$H^1_{\epsilon, mol} \hookrightarrow H^1_{L}.$$

\mb See also the work of P.Auscher, A.McIntosch and E. Russ \cite{AMR}, where a similar construction is done on a Riemannian manifold with doubling property, $L$ being the Hodge-De Rham Laplacian. One of the observation there is that one does not need pointwise bounds on the heat kernel but $L^2$ off-diagonal bounds similar to the ones we use in Section \ref{particular}. Again the Hardy space we obtain is included in the one of \cite{AMR}.

\gb {\bf Conclusion :} We have seen that our abstract construction is sometimes equal to sometimes smaller than other ones. We will see in Section \ref{secinter}
 that our space is big enough to interpolate with Lebesgue spaces. However, we think that
 this smallness is the main difficulty (which we will explain in Section \ref{secduality})
 for the identification of the dual spaces $(H^1_{\epsilon, mol})^{*}$ and $(H^1_{ato})^{*}$.

\section{Continuity theorem on the Hardy space.}

It is well known that a Calder\'on-Zygmund operator is continuous from the Coifman-Weiss
space $H^1_{CW}$ to $L^1$. We propose some general conditions which guarantee
the continuity from our Hardy spaces into $L^1$. We have the two following results~:

\begin{thm} \label{theo} Let $T$ be an $L^2$-bounded sublinear operator satisfying the following ``off-diagonal'' estimates~: for all ball $Q$, for all $j\geq 2$ there exist some coefficient $\alpha_j(Q)$ such that for all $L^2$-functions $f$ supported in $Q$
\be{hyp2h} \left(\frac{1}{\mu(2^{j+1}Q)} \int_{S_j(Q)} \left| T(B_Q(f))\right|^2 d\mu \right)^{1/2} \leq \alpha_j(Q) \left(\frac{1}{\mu(Q)}\int_{Q} |f|^2 d\mu \right)^{1/2}.\ee
If the coefficients $\alpha_j(Q)$ satisfy
\be{assum} \Lambda:=\sup_{Q \textrm{ ball}} \  \sum_{j\geq 2} \frac{\mu(2^{j+1}Q)}{\mu(Q)} \alpha_j(Q) <\infty, \ee
then there exists a constant $C$ such that
$$ \forall f\in L^2 \cap H^1_{ato} \qquad \|T(f)\|_{L^1} \leq C \|f\|_{H^1_{ato}}.$$
So if $T$ is linear then it has a unique extension, which is continuous from $H^1_{ato}$ into $L^1$.
\end{thm}

\begin{thm} \label{theo2h} Let $T$ be an $L^2$-bounded sublinear operator satisfying the following ``off-diagonal'' estimates~:
for all ball $Q$, for all $k\geq 0,j\geq 2$, there exist some coefficient $\alpha_ {j,k}(Q)$ such that for every $L^2$-function $f$ supported in $S_k(Q)$
 \be{assum1} \left(\frac{1}{\mu(2^{j+k+1}Q)} \int_{S_j(2^kQ)} \left| T(B_Q(f))\right|^2 d\mu \right)^{1/2} \leq \alpha_{j,k}(Q) \left(\frac{1}{\mu(2^{k+1}Q)}\int_{S_k(Q)} |f|^2 d\mu \right)^{1/2}. \ee
If the coefficients $\alpha_{j,k}$ satisfy \be{hyph} \Lambda:= \sup_{k\geq 0}\  \sup_{Q \textrm{ ball}} \
\left[\sum_{j\geq 2} \frac{\mu(2^{j+k+1}Q)}{\mu(2^{k+1}Q)} \alpha_{j,k}(Q) \right] <\infty, \ee then for all $\epsilon>0$ there exists a constant $C=C(\epsilon)$ such that
$$ \forall f\in L^2 \cap H^1_{\epsilon, mol} \qquad \|T(f)\|_{1} \leq C \|f\|_{H^1_{\epsilon, mol}}.$$
So if $T$ is linear then it has a unique extension, which is continuous from $H^1_{\epsilon, mol}$ to $L^1$.
\end{thm}

\begin{rem}
$1-)$ It is possible that a particular corona $S_k(Q)$ is empty. So we have normalized by $\mu(2^{k+1}Q)$. \\
$2-)$ One can weaken even more (\ref{hyp2h}) and (\ref{assum1}).
For example (\ref{hyp2h}) can be replaced by
$$ \int_{X \setminus 4Q} \left|T(B_Q(f))\right| d\mu \lesssim \mu(Q)^{1/2} \left(\int_Q |f|^2 d\mu \right)^{1/2}  $$
in which case the proof is almost tautological.
Also let us explain why we choose a condition like (\ref{hyp2h}).
Our Hardy spaces depend both on $\epsilon$ and on the collection
 ${\mathbb B} = (B_Q)_{Q \in {\mathcal Q}}$, so write them (just in this remark)
 $H^1_{\epsilon, mol, {\mathbb B}}$ and $H^1_{ato,{\mathbb B}}$.
 Take an $L^2$-bounded operator $T$. Assume that for $\epsilon'>0$ it
 satisfies the condition (\ref{hyp2h}) with some coefficients $\alpha_j(Q)$ satisfying~:
\be{assumee} \Lambda:=\sup_{Q \textrm{ ball}} \  \sum_{j\geq 2}
\frac{\mu(2^{j+1}Q)}{\mu(Q)} \alpha_j(Q)2^{j\epsilon'} <\infty. \ee
Then it is easy to show that for $m=B_Q(f_Q)$ an atom (of
$H^1_{ato,{\mathbb B}}$) associated to the ball $Q$, the function
$T(B_Q(f_Q))$ is an $\epsilon'$-molecule of $H^1_{\epsilon',T
{\mathbb B}}$ (associated to $Q$). Here we write
$T{\mathbb B}:= (TB_Q)_{Q\in{\mathcal Q}}$ the new collection of $L^2$-bounded operators. So we claim that $T$ is also continuous from $H^1_{ato,{\mathbb B}}$ into $H^1_{\epsilon',mol,T{\mathbb B}}$. We have an analogous result with (\ref{assum1}) and the molecular spaces. Also the assumptions (\ref{hyp2h}) and (\ref{assum1}) naturally appear when we want to work with these Hardy spaces. \\
$3-)$ Notice that when $\epsilon=\infty$, Theorem \ref{theo2h} becomes Theorem \ref{theo}. So it suffices to prove the last one.
\end{rem}

{\noindent {\bf Proof of Theorem \ref{theo2h}: }} First we show the following estimate : there exists a constant
$C=C(\epsilon)$ such that for all $\epsilon$-molecule $m$~:
\be{molborn} \left\|T(m) \right\|_{L^1} \leq C\left( \Lambda +
\|T\|_{L^2 \to L^2} \right). \ee
Using  definition we know that
there exists a ball $Q$ and a function $f_Q$ such that
$$m=B_Q(f_Q).$$
By decomposing the space $X$ with the scaled coronas around $Q$ and
by the linearity of $B_Q$, we have~:
$$ m=B_Qf_Q=\sum_{k\geq 0} B_Q({\bf 1}_{S_k(Q)} f_Q).$$
Using the sublinearity of $T$, we obtain that
$$\left|T(m)\right| \leq \sum_{k\geq 0} \left| TB_Q({\bf 1}_{S_k(Q)} f_Q) \right|.$$
By decomposing the integral with the coronas
$\left(S_j(2^kQ)\right)_{j\geq 0}$ which is a partition of $X$, we
have~:
\begin{align*}
\left\| T(m)\right\|_{L^1} & \leq \sum_{k\geq 0} \left\| TB_Q({\bf 1}_{S_k(Q)} f_Q) \right\|_{1} \\
  & \leq \sum_{\genfrac{}{}{0pt}{}{k\geq 0}{j\geq 0}} \int_{S_j(2^{k}Q)} \left| T(B_Q({\bf 1}_{S_k(Q)} f_Q))\right| d\mu \\
   & \leq \sum_{\genfrac{}{}{0pt}{}{k\geq 0}{j\geq 0}} \mu(2^{k+j+1}Q) \left(\frac{1}{\mu(2^{j+k+1}Q)} \int_{S_j(2^{k}Q)} \left| T(B_Q({\bf 1}_{S_k(Q)} f_Q))\right| d\mu \right) \\
   & \leq \sum_{\genfrac{}{}{0pt}{}{k\geq 0}{j\geq 0}} \mu(2^{k+j+1}Q) \left(\frac{1}{\mu(2^{j+k+1}Q)} \int_{S_j(2^{k}Q)} \left| T(B_Q({\bf 1}_{S_k(Q)} f_Q))\right|^2 d\mu \right)^{1/2}.
\end{align*}
Using the ``off-diagonal'' estimates (\ref{assum1}) on $T$ and the doubling condition for the measure $\mu$ (for the terms $j\leq 1$), we get
\begin{align*}
\left\| T(m)\right\|_{L^1} &  \leq \sum_{ \genfrac{}{}{0pt}{}{k\geq 0}{j\geq 2}} \mu(2^{k+1}Q) \frac{\mu(2^{k+j+1}Q)}{\mu(2^{k+1}Q)} \alpha_{j,k}(Q) \left(\frac{1}{\mu(2^{k+1}Q)} \int_{2^{k+1}Q} \left|{\bf 1}_{S_k(Q)} f_Q\right|^2 d\mu \right)^{1/2} \\
   & \qquad  + \sum_{k\geq 0,j\leq 1} A2^{j\delta} \mu(2^{k+1}Q) \left(\frac{1}{\mu(2^{k+1}Q)} \int_{X} \left|T(B_Q({\bf 1}_{S_k(Q)} f_Q))\right|^2 d\mu \right)^{1/2} \\
   & \leq \sum_{k\geq 0,j\geq 2}  \mu(2^{k+1}Q) \frac{\mu(2^{k+j+1}Q)}{\mu(2^{k+1}Q)} \alpha_{j,k}(Q) \mu(2^{k+1}Q)^{-1/2}  \left\|f_Q\right\|_{2,S_k(Q)} \\
   & \qquad  + \sum_{k\geq 0,j\leq 1}  A2^{j\delta} \mu(2^{k+1}Q)^{1/2}  \|TB_Q\|_{L^2\rightarrow L^2} \left\|f_Q\right\|_{2,S_k(Q)}.
\end{align*}
Then we use (\ref{operh}) to estimate $\|TB_Q\|_{L^2\rightarrow
L^2}$, and with the $L^2$-decay on $f_Q$, we have~:
\begin{align*}
 \left\| T(m)\right\|_{L^1}  & \leq \sum_{k\geq 0,j\geq 2}  \mu(2^{k+1}Q) \frac{\mu(2^{k+j+1}Q)}{\mu(2^{k+1}Q)} \alpha_{j,k}(Q) \mu(2^{k+1}Q)^{-1/2}  \mu(2^{k+1}Q)^{-1/2}2^{-\epsilon k} \\
     & \qquad \qquad \qquad  +  A'\|T\|_{L^2 \rightarrow L^2} \sum_{k\geq 0,j\leq 1}  2^{-\epsilon k+j\delta}  \\
   & \lesssim \sum_{k\geq 0} 2^{-\epsilon k} \left[\sum_{j\geq 2} \frac{\mu(2^{k+j+1}Q)}{\mu(2^{k+1}Q)} \alpha_{j,k}(Q) + 2^{\delta+1}\right]  \lesssim \Lambda + \|T\|_{L^2 \to L^2} .
\end{align*}
So we have proved the result for all $\epsilon$-molecules. \\
Next, we introduce the space
$$ S:= \left\{f;\  f=\sum_{i=1}^{n} \lambda_i m_i,\textrm{ $(m_i)_i$ $\epsilon$-molecules, } \|f\|_{H^1_{\epsilon, mol}} \geq 10^{-1} \left(\sum_{i=1}^n |\lambda_i|\right) \right\}, $$
which is a subspace of $H^1_{\epsilon, mol}$. \\
We have that $T$ is continuous from $S$ into $L^1$~: in fact if
$f\in S$, there exists a finite decomposition
$$ f = \sum_{i=1}^{n} \lambda_i m_i \qquad \textrm{with} \qquad \|f\|_{H^1_{\epsilon, mol}} \geq 10^{-1} \left(\sum_{i=1}^{n} |\lambda_i|\right), $$
by the sublinearity of $T$ and (\ref{molborn}), we have~:
$$\left\| T(f) \right\|_{L^1} \lesssim \sum_{i=1}^{n} |\lambda_i|\lesssim \|f\|_{H^1_{\epsilon,mol}}.$$
We conclude the proof by involving the next lemma. \findem

\begin{lem} \label{lemdense2} The set $S$ is dense in $H^1_{\epsilon, mol}$.
\end{lem}

\deme Let $f\in H^1_{\epsilon, mol}$ be a non zero function. Then there exists an infinite decomposition with $\epsilon$-molecules~:
$$f=\sum_{i=1}^\infty \lambda_i m_i \qquad \textrm{and} \qquad  \sum_{i=1}^\infty |\lambda_i| \leq (1+10^{-1})\|f\|_{H^1_{\epsilon, mol}}.$$
Let $f_{N}$ be the partial sums, defined by~:
$$f_N:=\sum_{i=1}^{N} \lambda_i m_i.$$
Assume that $f_N$ is not contained in $S$, then
$$\|f_{N} \|_{H^1_{\epsilon, mol}} \leq 10^{-1} \sum_{i=1}^{N} |\lambda_i|.$$
However $f_N$ converges to $f$ for the $H^1_{\epsilon, mol}$-norm, so for $N$ large enough we can deduce from this inequality that
$$\|f\|_{H^1_{\epsilon, mol}} \leq 2 \|f_{N} \|_{H^1_{\epsilon, mol}} \leq \frac{1}{5} \sum_{i=1}^{N} |\lambda_i| \leq \frac{1}{5} \sum_{i=1}^{\infty} |\lambda_i| \leq \frac{1}{2}\|f\|_{H^1_{\epsilon,mol}}.$$
This last inequality is not possible so we conclude that
for $N$ large enough, $f_N$ is an element of $S$. We have proved that $S$ is dense in $H^1_{\epsilon, mol}$. \findem

\begin{rem} The same lemma holds replacing $H^1_{\epsilon, mol}$ by $H^1_{ato}$.
\end{rem}

\begin{rem}  In order to verify the assumption (\ref{hyph}), with the doubling property of $\mu$ it is sufficient to check the stronger condition~:
\be{hyp3}  \sup_{k\geq 0} \  \sup_{Q \textrm{ball}} \  \left[\sum_{j\geq 2} 2^{j\delta} \alpha_{j,k}(Q) \right]
<\infty. \ee
\end{rem}

\begin{exmp} Theorem \ref{theo2h} (resp. Theorem \ref{theo}) applies to $T=Id$ if (\ref{assum1}) (resp. (\ref{hyp2h})) holds. In this case, (\ref{assum1}) becomes a condition on $B_Q$ which implies  \linebreak[4] $H^1_{\epsilon, mol} \hookrightarrow L^1$ (see Section \ref{particular}).
\end{exmp}

\begin{exmp} In the case of Coifman-Weiss Hardy space $H^1_{CW}$, this result generalizes the Calder\'on-Zygmund conditions. We choose
$$ B_Q(f)(x) = f(x) {\bf 1}_{Q}(x) - \left(\mu(Q)^{-1}\int_Q f d\mu\right) {\bf 1}_{Q}(x).$$
Let $T$ be an $L^2$ bounded operator. If we assume that the kernel $K(x,y)$ of $T$ satisfies the Calder\'on-Zygmund assumptions : there exists $h>0$ such that for all $x\neq y$ and $y'\in B(y,d(x,y)/2)$,
\begin{align*}
 &  \left| K(x,y) \right|  \lesssim \frac{1}{\mu(B(x,d(x,y)))}, & \\
 &  \left| K(x,y)-K(x,y')\right| \lesssim \frac{d(y,y')^h}{\mu(B(x,d(x,y))) d(x,y)^{h}}, &
\end{align*}
then it is easy to prove that $T$ satisfies the assumptions in
Theorem \ref{theo}. So $T$ is continuous from  $H^1_{CW}=H^1_{ato}$
into $L^1$, which is the classical result.
\end{exmp}

\mb So here we have obtained simple conditions for $L^2$-bounded
operators to be bounded from $H^1$ to $L^1$. These
conditions generalize the ``Calder\'on-Zygmund'' conditions in the
classical case. In next section, we are interested in interpolation
results.

\section{Interpolation theorem between $L^2$ and $H^1_{ato}$.} \label{secinter}

\mb The goal of this section is to find some general conditions on
$B_Q$ operators which give us an interpolation result like : if $T$,
an $L^2$-bounded operator, is continuous from $H^1$ to $L^1$ then
$T$ is $L^p$-bounded for all (or some) exponents $p\in]1,2]$. We
will use real $L^2$-duality defined by~:
$$ \forall f,g\in L^2, \qquad \langle f,g\rangle := \int f(x)g(x) d\mu(x).$$
Associated to this duality, we denote the adjoint operation by $^*$.

\begin{defn} We set $A_Q$ be the operator $Id-B_Q$ (here we can choose $A_Q(f)=f{\bf 1}_{Q} - B_Q(f)$ too). For $\sigma\in [1,\infty]$ we define the maximal operator~:
\be{opeM} \forall x\in X, \qquad M_{\sigma}(f)(x):= \sup_{\genfrac{}{}{0pt}{}{Q \textrm{ball}}{x\in Q}}\  \left( \frac{1}{\mu(Q)} \int_Q \left|A_Q^*(f)\right|^\sigma d\mu  \right) ^{1/\sigma} \ee
and a sharp maximal function adapted with our operators~: for $s>0$,
$$ \forall x\in X, \qquad M^\sharp_s (f)(x):= \sup_{\genfrac{}{}{0pt}{}{Q \textrm{ball}}{x\in Q}} \left(\frac{1}{\mu(Q)} \int_{Q} \left|B_Q^* (f)(z) \right|^{s} d\mu(z) \right)^{1/s}.$$
The standard maximal "Hardy-Littlewood" operator is defined by ~: for $s>0$,
$$\forall x\in X, \qquad M_{HL,s}(f)(x):=\sup_{\genfrac{}{}{0pt}{}{Q \textrm{ball}}{x\in Q}} \left(\frac{1}{\mu(Q)} \int_{Q} \left|f(z) \right|^{s} d\mu(z)
\right)^{1/s}.$$
\end{defn}

\mb For convenience, we recall the definition of a linearizable operator (see Definition V.1.20 of \cite{GR})~:
\begin{defn} \label{linearizable} An operator $T$ on $L^2$ is said to be linearizable if there exists a Banach space $\mathcal{B}$ and a linear operator $U$ defined from $L^2$ into $L^2(X,\mathcal{B})$ such that~:
$$\forall f\in L^2, \qquad \left| T(f)(x)\right| = \left\| U(f)(x)\right\|_{\mathcal{B}}, \qquad \mu-a.e.\ .$$
\end{defn}

\mb Examples of linearizable operators are given by maximal operators or
quadratic functionals for $T$. Linearizable operators are sublinear. Our main result is~:

\begin{thm} \label{theogeneh}  Let $\sigma\in]2,\infty]$. Assume that for all balls $Q$ and all functions $h\in L^2$, we have
$$ \left( \frac{1}{\mu(Q)} \int_Q |A_Q^*(h)|^\sigma d\mu \right)^{1/\sigma} \leq \inf_{x\in Q} M_{HL,2}(h)(x).$$
Let $T$ be an $L^2$-bounded, linearizable, operator. If $T$ is
continuous from $H^1_{ato}$ into $L^1$ then for all $p\in
]\sigma',2]$ there exists a constant $C=C(p)$ such that~:
$$ \forall f\in L^2 \cap L^p,\qquad \|T(f)\|_{p}\leq C\|f\|_{p}.$$
In addition we have the following estimate~:
$$\| T\|_{L^p \to L^p} \lesssim \|T\|_{L^2 \to L^2}^\theta \left[\|T\|_{H^1_{ato} \to L^1}^{1-\theta}+ \|T\|_{L^2 \to L^2}^{1-\theta}{\bf 1}_{\mu(X)<\infty} \right] $$
where $\theta$ is given by~:
$$\frac{1}{p} = \frac{\theta}{2}+ \frac{1-\theta}{1}.$$
The implicit constant depends on $A,A',\delta,p$ .
\end{thm}

\mb Here the quantity ${\bf 1}_{\mu(X)<\infty}$ means $1$ if $\mu(X)<\infty$ and $0$ if $\mu(X)=\infty$ and $\sigma'$ is the conjuguate exponent of $\sigma$. We need some preparation before proving this result.

\mb By Definition \ref{linearizable}, the operator $T$ is associated to a linear operator $U$ from $L^2$ into $L^2(X,\mathcal{B})$. We fix a measurable function $\phi$
$$ \phi : X \rightarrow \left\{ \lambda\in \mathcal{B}^{*},\ \|\lambda\|_{\mathcal{B}^*} \leq 1\right\}$$
and we consider the linear operator $V$ on $L^2$ defined by~:
\be{opV} \forall f\in L^2, \qquad V(f): \left\{ \begin{array}{lcl}
    X & \rightarrow & \R \\
    x & \mapsto & _{\mathcal{B}}\langle U(f)(x), \phi(x) \rangle_{\mathcal{B}^*}                                              \end{array} \right. . \ee
So $V$ is $L^2$-bounded  because
\begin{align*}
 \left\| V(f) \right\|_{2} & \leq \left\|_{\mathcal{B}}\langle U(f)(x), \phi(x) \rangle_{\mathcal{B}^*} \right\|_{2} \leq \left\| \left\|U(f)(x)\right\|_{\mathcal{B}} \left\|\phi(x) \right\|_{\mathcal{B}^*} \right\|_{2} \\
 & \leq \left\| \left\|U(f)(x)\right\|_{\mathcal{B}}\right\|_{2} = \left\| T(f) \right\|_{2}.
\end{align*}
By the same argument we can prove that $V$ is bounded from $H^1_{ato}$ to $L^1$ because we have~:
$$ \forall f\in H^1_{ato}, \qquad \left\| V(f) \right\|_{1} \leq \left\| T(f) \right\|_{1}.$$
By duality we know that $V^*$ is continuous from $L^\infty$ into $(H^1_{ato})^{*}$. So we have the following two lemmas~:
\begin{lem} \label{lem1} Let $C_1:=\|T\|_{H^1_{ato}\rightarrow L^1}$. Then ~:
$$\forall f\in L^\infty \cap L^{2} \qquad \left\| M^\sharp_2(V^*f) \right\|_\infty \leq C_1\|f\|_\infty.$$
\end{lem}

\deme Fix a function $f\in L^2 \cap L^\infty$. By the
$L^2$-boundedness of $V$ we have that $V^{*}(f)\in L^2$. Fix a ball
$Q$. Using the $L^2$-boundedness of $B_Q$ we conclude that
$B_Q^{*}(V^*f)$ exists and belongs to $L^2$. Let $h$ be supported
in $Q$ and normalized by $\|h\|_{2}=1$, and set $\phi_Q:=
\mu(Q)^{-1/2}h$. Then it is easy to see that $m=B_Q(\phi_Q)$ is an
atom. With the continuity of $V$ from $H^1_{ato}$ to $L^1$, we
obtain~:
$$\left| \int V(m) f d\mu \right| \leq \|V\|_{H^1_{ato}\rightarrow L^1} \|f\|_\infty\leq C_1 \|f\|_\infty.$$
In addition, since $m$ and $f$ are $L^2$-functions, we have~:
\begin{align*}
\int V(m) f d\mu & = \int m V^{*}(f) d\mu = \int \mu(Q)^{-1/2} B_Q(h) V^{*}(f) d\mu \\
 & = \mu(Q)^{-1/2} \int h B_Q^{*} V^* (f)d\mu .
\end{align*}
In consequence, we get
$$ \forall h\in L^2(Q),\ \|h\|_{2,Q}=1, \qquad \left| \int h B_Q^* V^*(f) d\mu \right| \leq C_1 \mu(Q)^{1/2} \|f\|_\infty.$$
Hence, $B_Q^*V^*(f) \in L^2(Q)$ and
$$\left\| B_Q^* V^*(f) \right\|_{2,Q} \leq C_1 \mu(Q)^{1/2} \|f\|_\infty,$$
which concludes the proof. \findem

\begin{lem} \label{lem2} For all $1\leq s\leq 2$, there exists a constant $C_2=C_2(s)\lesssim \|T\|_{L^2 \to L^2}$ such that~:
$$\forall f\in L^{2} \qquad \left\| M^\sharp_s(V^*f) \right\|_{2,\infty} \leq C_2\|f\|_2.$$
\end{lem}

\deme Let $x\in X$ and $Q$ be a ball such that $x\in Q$. We have~:
\begin{align*}
\lefteqn{\left(\frac{1}{\mu(Q)} \int_{Q} \left|B_Q^* (V^* f)(z) \right|^{s} d\mu(z) \right)^{1/s}} & & \\
 & \qquad  \leq \left(\frac{1}{\mu(Q)} \int_{Q} \left| (V^*f)(z) \right|^{s} d\mu(z)\right)^{1/s}  + \left(\frac{1}{\mu(Q)} \int_{Q} \left|A_Q^* (V^*f)(z) \right|^{\sigma} d\mu(z)\right)^{1/\sigma} & \\
 & \qquad \qquad \leq M_{HL,s}(V^*f)(x) + M_\sigma(V^*f)(x). &
 \end{align*}
Here we use the fact that $s\leq 2\leq \sigma$. Taking the supremum
over all balls $Q \ni x$, we get~:
$$M^\sharp_s(V^*f) \leq M_{HL,s}(V^*f) + M_\sigma(V^*f).$$
In addition for $s\leq 2$, $M_{HL,s}$ is of weak type $(2,2)$.
 By the assumptions $V^*$ is $L^2$-bounded and $M_\sigma$ is of weak type $(2,2)$ (bounded by $M_{HL,2}$). So we conclude that $M^\sharp_s \circ V^*$ is of weak type $(2,2)$.
\findem

\mb Let the parameter $s=2$. From the two previous lemmas (Lemma
\ref{lem1} and \ref{lem2}) and Marcinkiewicz's theorem for real
interpolation, we obtain that $M^\sharp_2 \circ V^*$ is bounded in
the whole space $L^r$ for $2< r < \infty$. If $p\in]1,2[$ is an
exponent and $q\in]2,\infty[$ its conjuguate exponent, then there
exists also a constant $C_3:=C_3(p,C_1,C_2)$ such that \be{conth}
\forall f\in L^{q}\cap L^2 \qquad \left\| M^\sharp_2(V^*f)
\right\|_{q} \leq C_3 \|f\|_{q}, \ee and $C_3$ is bounded by~:

\begin{align*}
C_3 & \lesssim \left\|M^\sharp_2 \circ V^* \right\|_{L^2 \to L^{2,\infty}}^{\theta} \left\|M^\sharp_2 \circ V^*
\right\|_{(H^1_{ato})^* \to L^\infty}^{1-\theta} \\
 & \lesssim \left\|T \right\|_{L^2 \to L^{2}}^{\theta} \left\|T
\right\|_{H^1_{ato} \to L^1}^{1-\theta},
\end{align*}
because $\theta$ satisfies
$$\frac{1}{q}= 1-\frac{1}{p}=\frac{\theta}{2}.$$

\mb We now begin the proof of Theorem \ref{theogeneh}.

\deme We develop the proof in two steps. \\
$1-)$ End of the proof for the operator $V$. \\
We are going to obtain the result as a consequence of Theorem 3.1 in \cite{AM}. With its notations, take $F=|h|^2$ and for all balls $Q$
$$ G_Q= 2|B_Q^*h|^2 \textrm{    and    } H_Q= 2|A_Q^*h|^2,$$
where $h=V^{*}(f)$ and $f\in L^2$.
Then for all $x\in Q$, we have $F(x)\leq G_Q(x)+H_Q(x)$. For all balls $Q \ni x$, we have
$$ \frac{1}{\mu(Q)} \int_Q G_Q d\mu \leq 2 M_{2}^\sharp(h)(x)^2:=G(x).$$
By assumption, for all balls $Q \ni x$ and $\overline{x}\in Q$, we have
\begin{align*} \left( \frac{1}{\mu(Q)} \int_{Q} |H_Q|^{\sigma/2} d\mu \right)^{2/\sigma} & \leq 2^{2/\sigma} \left( \frac{1}{\mu(Q)} \int_Q |A_Q^*(h)|^\sigma d\mu \right)^{2/\sigma} \leq 2^{2/\sigma} M_{\sigma}(h)(\overline{x})^2 \\
& \leq 2^{2/\sigma} M_{HL,2}(h)(\overline{x})^2 = 2^{2/\sigma} M_{HL,1}(F)(\overline{x}).
\end{align*}
Also by Theorem 3.1 of \cite{AM} (which applies because $h=V^*(f) \in L^2$ so $F\in L^1$), for $1<r<\sigma/2$ we have~:
\be{cont11} \|F\|_{r} \leq \| M_{HL,1}(F) \|_{r} \lesssim \| M_2^\sharp(h)^2 \|_{r} = \|M_2^\sharp(h) \|_{2r}^2. \ee
By using (\ref{conth}) and $h= V^*(f)$, we obtain for $2 < 2r<\infty$~:
\be{cont110} \| M_2^\sharp(h) \|_{2r} \lesssim \|f\|_{2r}. \ee
Then the two previous estimates (\ref{cont11}) and (\ref{cont110}) give us
$$ \| V^{*}(f) \|_{2r} \lesssim \|f\|_{2r}.$$
We have obtained also that for all exponents $q$ with $2<q<\sigma$, there exists a constant $C=C(q)$ such that~:
$$\forall f\in L^2 \cap L^q, \qquad \|V^{*}(f) \|_{q} \leq C \|f\|_{q}.$$
Thus by duality we can obtain that $V$ is bounded on $L^p$ for all exponents $p\in]\sigma',2[$. \\
$2-)$ End of the proof for the operator $T$. \\
We have obtained that there exists a constant $C_4$ such that for all functions $\phi$ and for all functions $f\in L^p\cap L^2$ we have \be{contV}
\left\| V(f) \right\|_{p} \leq C_4 \|f\|_{p}. \ee Using the Hahn-Banach theorem, we have that
$$\forall b\in \mathcal{B}, \exists b^*\in \mathcal{B}^*  \qquad \frac{1}{4} \left\|b \right\|_{\mathcal{B}} \leq  _{\mathcal{B}} \langle b, b^* \rangle_{\mathcal{B}^*} \qquad \textrm{ and } \qquad \|b^{*}\|_{\mathcal{B}^*} \leq 1. $$
By choosing a good measurable function $\phi$ (see Remark \ref{remden}), we can have also
 $$  V(f)(x) \geq \frac{1}{4} \left\| U(f)(x) \right\|_{\mathcal{B}} = \frac{1}{4} \left| T(f)(x)\right|.$$
Hence by (\ref{contV}), we obtained that
$$ \left\| T(f) \right\|_{p} \lesssim \|f\|_{p}.$$
This estimate is uniform for all $f\in L^p \cap L^2$.
\findem

\begin{rem} \label{remden} By density we can approach the $\mathcal{B}$-valued function $U(f)(x)$ by a finite sum $\sum {\bf 1}_{A_i}(x) b_i$ with $A_i$ mesurable sets of $X$ and $(b_i)_i$ a finite independant family of ${\mathcal B}$. Then we chose the function $\phi(x):=\sum {\bf 1}_{A_i} b_i^*$ where $b_i^*\in {\mathcal B}^{*}$ are chosen so that $\langle b_j,b_i^{*} \rangle = \|b_i\|_{\mathcal B}{\bf 1}_{i=j}$, in order to have $$ \left\| \|U(f)(x)\|_{{\mathcal B}} \right\|_{p,d\mu(x)} \simeq \left\| \langle T(f)(x), \phi \rangle \right\|_{p,d\mu(x)}.$$
\end{rem}

\mb There is an interesting particular case when $T=Id$, which gives us the following corollary~:

\begin{cor} \label{corint} Assume that $M_{\sigma}$ is bounded by $M_{HL,2}$ for some $\sigma\in]2,\infty]$ and that  $H^1_{ato}$ is continuously contained in $L^1$ (see Section \ref{particular}). Then for all $q\in ]2,\infty[$ there exists a constant $C=C(q)$ such that
$$ \forall f\in L^q\cap L^2, \qquad  \| M^\sharp_2(f) \|_{q} \leq C \|f\|_{q}.$$
In addition if $1\leq s \leq 2< q <\sigma$ then we have the other inequality~:
$$ \forall f\in L^q\cap L^2, \qquad C^{-1} \|f\|_{q} \leq \| M^\sharp_s(f) \|_{q} \leq C\|f\|_{q}.$$
\end{cor}

\deme With these assumptions, we can apply the previous theorem with
$T=Id$. The first claim is obvious with (\ref{conth}), here we have
just used the imbedding $H^1_{ato} \hookrightarrow L^1$. The second
claim is a consequence of (\ref{cont11}) with the control $|f|\leq
M_{HL,s}(f)$. To prove this, we use the assumption that $M_\sigma$
is bounded by $M_{HL,2}$. In the proof of Theorem \ref{theogeneh}, we have used $s=2$
however we can use the same arguments with $s\in[1,2]$. \findem

\mb We have a second corollary (independent of Corollary \ref{corint})~:

\begin{cor} \label{corinter}  Assume that for $\sigma\in ]2,\infty]$ the maximal operator $M_{\sigma}$ is bounded by $M_{HL,2}$. Then we have the following
inequality : for all $q\in ]2,\sigma[$ there exists a constant $c_q$ such that for all functions $f\in L^2 \cap
(H^1_{ato})^*$ we have \be{coreq} \| f\|_{q} \leq c_q \|f\|_{2}^{\theta} \left[ \|f\|_{(H^1_{ato})^*}^{1-\theta}+
\|f\|_{2}^{1-\theta}{\bf 1}_{\mu(X)<\infty} \right], \ee where $\theta$ is given by
$$\frac{1}{q} = \frac{\theta}{2}.$$
\end{cor}

\deme In fact we can prove this result  directly by using the
maximal operator $M^\sharp_2$ and the previous arguments. Here we
will prove this as an application of the previous theorem. So we
take a subset $E$ of $X$ satisfying $0<\mu(E)<\infty$, and we write
$\phi:={\bf 1}_{E}$. We define also the operator $T$ as~:
$$T(h) := \langle h , f\rangle \phi.$$
The assumption ``$f\in L^2$'' guarantees that $T$ is $L^2$-bounded. And by the fact that $f\in (H^1_{ato})^*$, we obtain that $T$ is continuous from $H^1_{ato}$ to $L^1$. So we can apply the previous theorem and we conclude
that $T$ is $L^p$ bounded for all $p\in ]\sigma',2[$, and we have the following estimate~:
$$\left\| T(h) \right\|_{p} \lesssim \|h\|_{p} \| T\|_{L^2 \to L^2}^{\theta} \left[\|T\|_{H^1_{ato} \to L^1}^{1-\theta}+ \| T\|_{L^2 \to L^2}^{1-\theta}{\bf 1}_{\mu(X)<\infty}\right].$$
By the definition of $T$ and the duality result about Lebesgue spaces, we deduce that (with $p$ equals to the conjuguate exponent of $q$)~:
$$\|f\|_{q} \simeq \|T\|_{L^p \to L^p} \lesssim \|f\|_{2}^\theta \left[\|f\|_{(H^1_{ato})^*}^{1-\theta}+\|f\|_{2}^{1-\theta}{\bf 1}_{\mu(X)<\infty}\right].$$
\findem

\mb The parameter $\sigma$ plays an important role in Theorem \ref{theogeneh}. It permits to understand which Lebesgue spaces $L^p$ can be obtained by interpolation between the Hardy space $H^1_{ato}$ and $L^2$. We now see that the range for $p$ is optimal with the following example.

\begin{exmp} \label{exhm} The Riesz transforms. \\
Here take $X=\R^n$ and let
$$ A : \R^n \longrightarrow {\mathcal M}_{n}(\C) $$
 be a bounded $n \times n$ matrix valued function satisfying the ellipticity condition : there exist two constants $\Lambda\geq \lambda>0$ such that
$$ \forall x\in \R^n, \ \forall \xi,\zeta\in \C^n, \qquad \lambda|\xi|^2 \leq Re \left(A(x)\xi \cdot \overline{\xi} \right) \quad  \textrm{and} \quad |A(x)\xi \cdot \overline{\zeta} | \leq \Lambda|\xi||\zeta|.$$ We define the second order divergence form operator as
$$ L(f):= -\textrm{div} (A \nabla f).$$
Associated to this operator we have the Riesz transform $\nabla L^{-1/2}$ defined by
$$\nabla L^{-1/2}(f) := \frac{1}{\sqrt{\pi}} \int_0^\infty \sqrt{t} \nabla e^{-tL} (f) \frac{dt}{t}.$$
From \cite{A}, we know that the interval of $p\in[1,2]$ such that the heat semigroup $(e^{tL})_{t\geq 0}$ and the Riesz transform  $\nabla L^{-1/2}$ are $L^{p}$-bounded have the same interior. We shall denote the critical exponent $p_{L}$. Moreover we know that for $p\in[1,2]$
$$ \nabla L^{-1/2} \textrm{ is $L^p$-bounded} \Longleftrightarrow p_{L}<p\leq 2.$$
It is shown in \cite{A1} that $1\leq p_{L}<\frac{2n}{n+2}$ and in
\cite{A} that it could be that $p_L>1$. In that case the Riesz transform
can not be continuous from $H^1_{CW}$ into $L^1$ (otherwise
interpolation would yield $p_L=1$). So this is an example where the
Coifman-Weiss Hardy space is not well adapted to the operator.

\mb We have seen in Subsection \ref{elliptic} that we can compare the Hardy space $H^1_{L}$ of \cite{HM} with ours. By choosing
$$ B_Q(f) := (r_Q^2L)^{M}e^{-r_Q^2L}(f) \textrm{  or  }  B_Q(f) := \left(Id - (Id+r_Q^2L)^{-1} \right)^{M}(f),$$
with $M$ a large enough integer ($M>n/4$), we have the inclusion
\be{exinclus} \forall \epsilon>0, \qquad H^1_{\epsilon, mol} \hookrightarrow H^1_{L}. \ee \\
A "weak" version of Lemma 2.5 in \cite{HM} asserts that for every $q$, $2\leq q < (p_L)'$, for arbitrary closed sets $E,F\subset \R^n$
$$ \forall t>0, \forall f\in L^2(E), \qquad \|(tL^*)^{M}e^{-tL^*}(f) \|_{q,F} \lesssim t^{\frac{n}{2}(\frac{1}{q}-\frac{1}{2})} e^{-cd(E,F)^2/t} \|f\|_{2,E}.$$

\mb By taking $t=r_Q^2$, $q=\sigma<(p_L)'$, it is easy to check that for all balls $Q$, for all $h\in L^2$
$$ \left( \frac{1}{|Q|} \int_Q \left| (r_Q^2L^*)^{M}e^{-r_Q^2L^*}(h)(x) \right|^\sigma dx \right)^{1/\sigma} \lesssim \inf_{z\in Q} M_{HL,2}(h)(z).$$
This gives us that the maximal operator $M_{\sigma}$ (defined by
(\ref{opeM})) satisfies~:
$$\forall h\in L^2, \qquad M_\sigma(h) \lesssim M_{HL,2}(h).$$
We already know that $\nabla L^{-1/2}$ is $L^2$-bounded (\cite{AHLMT}). For $\epsilon>0$,
Theorem 3.4 of \cite{HM} and (\ref{exinclus}) prove that this operator is continuous
from $H^1_{\epsilon, mol}$ into $L^1$. We can also apply our Theorem \ref{theogeneh}
to deduce that the Riesz transform is $L^p$-bounded for all $p\in]\sigma',2]$
with all $\sigma\in[2,p_L'[$. So the Riesz transform is $L^p$ bounded for
all exponents $p\in ]p_L,2]$. As we know the Riesz transform is not $L^p$-bounded
for $p\leq p_L$ so our range of exponents is optimal.
\end{exmp}

\mb From the point of view of "$L^p$-theory" (we temporary forget the Hardy spaces),
we have the following result (using Theorem \ref{theo2h} and Theorem \ref{theogeneh})~:

\begin{thm} \label{theo4} Suppose $1\leq p_0<2$. Let $T$ be an $L^2$-bounded, linearizable, operator such that for all balls $Q$ and for all functions $f$ supported in $Q$
$$ \forall j\geq 2 \quad \left(\frac{1}{\mu(2^{j+1}Q)} \int_{S_j(Q)} \left| T(B_Q(f))\right|^2 d\mu \right)^{1/2} \leq \alpha_j(Q) \left(\frac{1}{\mu(Q)}\int_{Q} |f|^2  d\mu \right)^{1/2}$$
and
$$\forall j\geq 0 \quad \left(\frac{1}{\mu(2^{j+1}Q)} \int_{S_j(Q)} \left| f-B_Q(f)\right|^2 d\mu \right)^{1/2} \leq \alpha_j(Q) \left(\frac{1}{\mu(Q)}\int_{Q} |f|^{p_0} d\mu \right)^{1/p_0},$$
where the coefficients $\alpha_j(Q)$ satisfy
$$\sup_{Q \textrm{ ball}}\ \sum_{j\geq 0} \frac{\mu(2^{j+1}Q)}{\mu(Q)} \alpha_j(Q) <\infty.$$
Then for all exponents $p\in ]p_0,2[$, there exists a constant $C$ such that
$$ \forall f\in L^2 \cap L^p, \qquad  \|T(f)\|_{p} \leq C \|f\|_{p}.$$
\end{thm}

\deme We are going to show that with these assumptions the maximal
operator $M_{(p_0)'}$ (defined by (\ref{opeM})) is bounded by
$M_{HL,2}$. Then the theorem is a consequence of Theorem
\ref{theo2h} and Theorem \ref{theogeneh}. To prove this, let $f\in
L^2\cap L^{1}$ be a function and $\sigma=(p_0)'$.
 For all balls $Q$ containing the point $x_0$, we can estimate the $L^{\sigma}$-norm by duality~:
\begin{align*}
M_\sigma (f)(x_0)& =\sup_{x_0 \in Q} \sup_{\genfrac{}{}{0pt}{}{g \in L^2(Q)}{\|g\|_{p_0}\leq 1}} \mu(Q)^{-1/\sigma} \int \left( f - B_Q^*(f)\right) g d\mu \\
 & = \sup_{x_0 \in Q} \sup_{\genfrac{}{}{0pt}{}{g \in L^2(Q)}{\|g\|_{p_0}\leq 1}} \mu(Q)^{-1/\sigma} \int f \left( g - B_Q(g) \right) d\mu \\
 & \leq \sup_{x_0 \in Q} \sup_{\genfrac{}{}{0pt}{}{g \in L^2(Q)}{\|g\|_{p_0}\leq 1}} \mu(Q)^{-1/\sigma} \sum_{j\geq 0} \int_{S_j(Q)} \left|f \left(g-B_Q(g) \right) \right| d\mu  \\
 & \leq \sup_{x_0 \in Q} \mu(Q)^{-1/\sigma} \sup_{\genfrac{}{}{0pt}{}{g \in L^2(Q)}{\|g\|_{p_0}\leq 1}} \sum_{j\geq 0}  \left\|f \right\|_{2,S_j(Q)} \left\| g-B_Q(g) \right\|_{2,S_j(Q)} \\
 & \leq \sup_{x_0 \in Q} \mu(Q)^{-1/\sigma} \sup_{\genfrac{}{}{0pt}{}{g \in L^2(Q)}{\|g\|_{p_0}\leq 1}}  \sum_{j\geq 0}  \left\|f \right\|_{2,S_j(Q)} \alpha_j(Q) \mu\left(2^{j+1}Q\right)^{1/2} \mu(Q)^{-1/p_0}\left\|g\right\|_{p_0} \\
& \leq \sup_{x_0 \in Q} \sum_{j\geq 0}  \left\|f \right\|_{2,S_j(Q)} \mu(2^{j+1} Q)^{-1/2} \alpha_j(Q) \frac{\mu(2^{j+1}Q)}{\mu(Q)} \\
& \leq M_{HL,2}(f)(x_0) \sup_{x_0 \in Q}  \sum_{j\geq 0}   \alpha_j(Q) \frac{\mu(2^{j+1}Q)}{\mu(Q)} \\
& \lesssim M_{HL,2}(f)(x_0).
\end{align*}
\findem

\mb Let us compare with S.Blunck and P.Kunstmann's result (see
Theorem 1.1 of \cite{bk}). We describe an improved version, due to P.
Auscher (Theorem 1.1 of \cite{A})~:

\begin{thm} \label{theo5} Let $1\leq p_0 <2$. Let $\widetilde{B}_Q$ be some uniformly $L^2(X)$-bounded operators. Let $T$ be a $L^2$-bounded sublinear operator such that for all balls $Q$ and for all functions $f$ supported in $Q$
$$\forall j\geq 2 \quad\left(\frac{1}{\mu(2^{j+1}Q)} \int_{S_j(Q)} \left| T(\widetilde{B}_Q(f))\right|^2 d\mu \right)^{1/2} \leq \alpha_j(Q) \left(\frac{1}{\mu(Q)}\int_{Q} |f| d\mu \right)$$
and
$$ \forall j\geq 0 \quad \left(\frac{1}{\mu(2^{j+1}Q)} \int_{S_j(Q)} \left| f-\widetilde{B}_Q (f))\right|^2 d\mu \right)^{1/2} \leq \alpha_j(Q) \left(\frac{1}{\mu(Q)}\int_{Q} |f|^{p_0} d\mu \right)^{1/p_0}.$$
If the coefficients $\alpha_j(Q)$ satisfy
$$\sup_{Q \textrm{ ball}} \sum_{j\geq 0} \frac{\mu(2^{j+1}Q)}{\mu(Q)} \alpha_j(Q) <\infty,$$
then there exists a constant $C$ such that
$$ \forall f\in L^2 \cap L^{p_0}, \qquad  \|T(f)\|_{L^{p_0,\infty}} \leq C \|f\|_{p_0}.$$
So by real interpolation, for all $p\in]p_0,2[$, there exists a constant $C$ such that
$$ \forall f\in L^2 \cap L^p, \qquad  \|T(f)\|_{p} \leq C \|f\|_{p}.$$
\end{thm}

\begin{rem} \label{remkb} The first ``$L^{1}-L^2$'' condition of Theorem \ref{theo5}
is stronger than the first ``$L^2-L^2$'' condition of Theorem
\ref{theo4}. So Theorem \ref{theo4} has weaker conditions as far as
the continuity on $L^p$ with any $p\in ]p_0,2[$ is concerned.
However the assumptions in Theorem \ref{theo5} imply the weak type
$(p_0,p_0)$ by a variant of the Calder\'on-Zygmund decomposition,
which does not seem to work under the assumptions of Theorem
\ref{theo4}. Also Theorem \ref{theo5} can be applied to sublinear
operators, not only to linearizable operators.
\end{rem}

\section{Continuity results in weighted spaces.}

We have used the sharp maximal function $M_s^\sharp$ in previous
sections. Like in the Euclidean case, we can use it to get
weighted norm inequalities.

\mb We have the following result~:

\begin{thm} \label{weight} Let $T$ be a linear operator as in Theorem \ref{theo}, then for $1\leq s\leq 2$, we have for $f\in L^2$,
$$ M^\sharp_{s}(T^*(f)) \lesssim M_{HL,2}(f).$$
\end{thm}

\deme Let $f$ be a function in $L^2$. Let us fix $x_0$ a point and write $s'$ the conjuguate exponent of $s$. By using the assumption (\ref{hyp2h}) of Theorem \ref{theo}~:
\begin{align*}
M^\sharp_s(T^*f)(x_0) & =\sup_{x_0 \in Q} \sup_{\genfrac{}{}{0pt}{}{g \in L^2(Q)}{\|g\|_{s'}\leq 1}} |Q|^{-1/s}\int B_Q^*(T^*f) g d\mu \\
 & = \sup_{x_0 \in Q} \sup_{\genfrac{}{}{0pt}{}{g \in L^2(Q)}{\|g\|_{s'}\leq 1}} |Q|^{-1/s}\int f T(B_Q(g)) d\mu \\
 & \leq \sup_{x_0 \in Q} \sup_{\genfrac{}{}{0pt}{}{g \in L^2(Q)}{\|g\|_{s'}\leq 1}} \sum_{i\geq 0} \mu(Q)^{-1/s}\int_{S_i(Q)} \left|f T(B_Q(g)) \right| \\
 & \leq \sup_{x_0 \in Q} \sup_{\genfrac{}{}{0pt}{}{g \in L^2(Q)}{\|g\|_{s'}\leq 1}} \sum_{i\geq 0}  \mu(Q)^{-1/s}\left\|f \right\|_{2,S_i(Q)} \left\| T(B_Q(g))
 \right\|_{2,S_i(Q)}.
\end{align*}
For $x_0\in Q$, we have the following estimates~:
\begin{align*}
 \sum_{i\geq 0} & \mu(Q)^{-1/s}\left\|f \right\|_{2,S_i(Q)} \left\| T(B_Q(g))\right\|_{2,S_i(Q)}  \\
 &  \leq  \mu(Q)^{-1/s}\left[ \sum_{i \leq 2}\left\|f \right\|_{2,S_i(Q)} \|g\|_2 + \sum_{i\geq 3}  \left\|f \right\|_{2,S_i(Q)} \alpha_i(Q) \left(\frac{\mu \left(2^{i+1}Q\right)}{\mu(Q)} \right)^{1/2}\left\|g\right\|_{2} \right] \\
&  \lesssim  \mu(Q)^{1/2}\|g\|_{s'}\left[ \sum_{i \leq 2}\left\|f \right\|_{2,S_i(Q)}   + \sum_{i\geq 3}  \left\|f \right\|_{2,S_i(Q)} \alpha_i(Q) \left( \frac{\mu \left(2^{i+1}Q\right)}{ \mu(Q)} \right)^{1/2}  \right] \\
 &  \lesssim M_{HL,2}(f)(x_0)  \left[1 + \sum_{i \geq 3}  \alpha_i(Q) \mu \left(2^{i+1}Q\right) \mu(Q)^{-1} \right] \\
&  \lesssim M_{HL,2}(f)(x_0).
\end{align*}
So we obtain~: $$ M^\sharp_s(T^*f)(x_0) \lesssim M_{HL,2}(f)(x_0).$$
 \findem

\mb We recall the definition of Muckenhoupt's weights and Reverse Hölder classes~:

\begin{defn} A nonnegative function $\omega$ on $X$ belongs to the class ${\mathbb A}_p$ for $1< p <\infty$ if
$$\sup_{Q \textrm{ ball}} \left(\frac{1}{\mu(Q)}\int_Q w d\mu \right) \left( \frac{1}{\mu(Q)}\int_Q \omega^{-1/(p-1)} d\mu \right)^{p-1} <\infty.$$
A nonnegative function $\omega$ on $X$ belongs to the class $RH_q$ for $1<q<\infty$, if there is a constant $C$ such that for every ball $Q\subset X$
$$ \left( \frac{1}{\mu(Q)}\int_Q \omega^q d\mu \right)^{1/q} \leq C \left(\frac{1}{\mu(Q)}\int_Q \omega d\mu \right).$$
\end{defn}

\mb We have the well-known following properties (chapter 9 of \cite{Gra2} for the Euclidean case)~:

\begin{prop} \label{prora} For $1<s<\infty$ the maximal operator $M_{HL,s}$ is bounded on $L^p(\omega)$ for all $s<p<\infty$ and $\omega\in {\mathbb A}_{p/s}$. For $1\leq p<\infty$ and $\omega$ an ${\mathbb A}_p$-weight, there exists some constants $C,\epsilon>0$ such that for all balls $Q$ and all mesurable subsets $A\subset Q$, we have:
 \be{poidsh} \frac{\omega(A)}{\omega(Q)} \leq C \left(\frac{\mu(A)}{\mu(Q)} \right)^{\epsilon}. \ee
If $\omega$ is as ${\mathbb{A}}_p$ weight and $1<p<\infty$, then $\omega^{1-p'}=\omega^{-1/(p-1)}$ is an ${\mathbb{A}}_{p'}$ weight. In addition we have the following equivalence~:
$$ \omega \in{\mathbb{A}}_{p} \Longleftrightarrow \omega^{1-p'} \in {\mathbb{A}}_{p'}.$$
\end{prop}

\mb With these definitions we have the following result~:

\begin{thm} \label{poidsthm} Let $T$ be a linearizable operator which satisfies the assumption of Theorem \ref{theo}. We assume that for $\sigma\in ]2,\infty]$ the maximal operator $M_\sigma$ is bounded by $M_{HL,2}$. Let $p\in]\sigma',2[$ be an exponent and $\omega$ a weight so that $\omega \in {\mathbb A}_{p/\sigma'} \cap RH_{(\frac{2}{p})'}$. Then the operator $T$ is "bounded" on $L^p(\omega)$ : there exists a constant $C$ such that~:
$$ \forall f\in L^2 \cap L^p(\omega), \qquad \|T(f)\|_{p,\omega d\mu} \leq C\|f\|_{p,\omega d\mu}.$$
\end{thm}

\mb For the weight $\omega$, we define the associated measure (written by the same symbol) $\omega$ by $d\omega:=\omega d\mu$. Concerning the condition $\omega \in {\mathbb A}_{p/2} \cap RH_{(\frac{\sigma'}{p})'}$, we recall (Lemma 4.4 of \cite{AM}) that for $\sigma'<p<2$
\be{poidseq}
\omega \in {\mathbb A}_{p/\sigma'} \cap RH_{(\frac{2}{p})'} \Longleftrightarrow \omega^{1-p'}\in {\mathbb A}_{p'/2}\cap RH_{(\frac{\sigma}{p'})'}.
\ee
In the following proof, we will use this equivalence.

\mb
\deme We use the notation of Theorem \ref{theogeneh}. Let $f\in L^2 \cap L^{q}(\omega^{1-q})$ where $q=p'$ is the conjuguate exponent of $p$. With $h=V^{*}(f)$, we write
$F=|h|^2$. \\
First we recall the fact that $V^*$ is the adjoint of $V$ related to the measure $\mu$. So we have that
$$ \textrm{ $V$ is $L^p(\omega)$-bounded} \ \Longleftrightarrow \ \textrm{ $V^*$ is $L^q(\omega^{1-q})$-bounded}.$$
We use Theorem 3.1 of \cite{AM} with $r:=q/2 \in ]1,\sigma/2[$ and $\omega^{1-2r} = \omega^{1-p'} \in RH_{(\frac{\sigma}{p'})'}$ to obtain the weighted version of (\ref{cont11})~:
\be{cont11w} \|F\|_{r,\omega^{1-2r}} \leq \| M_{HL,1}(F) \|_{r,\omega^{1-2r}} \lesssim \| M_2^\sharp(h)^2 \|_{q,\omega^{1-2r}} = \|M_2^\sharp(h) \|_{2r,\omega^{1-2r}}^2. \ee
By using Theorem \ref{weight} and $h= V^*(f)$, we obtain for $2 < q=2r <\infty$~:
$$\| M_2^\sharp(h) \|_{2r,\omega^{1-2r}} \lesssim \|M_{HL,2}(f)\|_{2r,\omega^{1-2r}}.$$
The weight $\omega^{1-q}$ belongs to the class ${\mathbb A}_{q/2}$, so with Proposition \ref{prora} we get
$$ \|M_{HL,2}(f)\|_{q,\omega^{1-q}} \lesssim \|f\|_{q,\omega^{1-q}}.$$
Then the three previous estimates give us
$$ \| V^{*}(f) \|_{q,\omega^{1-q}} \lesssim \|f\|_{q,\omega^{1-q}}.$$
Thus by duality we obtain that $V$ is bounded on $L^p(\omega)$ and we deduce the same for $T$ as in the proof of Theorem \ref{theogeneh}.
\findem

\mb We use the following notation of \cite{AM}~:

\begin{defn} For $\omega$ a nonnegative function on $X$ and $0<p_0<q_0\leq \infty$ two exponents, we introduce the set
$$ \mathcal{W}_{\omega}(p_0,q_0):=\left\{ p\in]p_0,q_0[,\ \omega\in {\mathbb A}_{p/p_0} \cap RH_{(q_0/p)'} \right\}.$$
\end{defn}

\mb Then we have the "weighted" version of Theorem \ref{theo4}~:

\begin{thm} \label{theo4poids} Let $1\leq p_0 <2$. Let $T$ be an $L^2$-bounded, linearizable, operator such that for all balls $Q$ and for all functions $f$ supported in $Q$
$$ \forall j\geq 2 \quad \left(\frac{1}{\mu(2^{j+1}Q)} \int_{S_j(Q)} \left| T(B_Q(f))\right|^2 d\mu \right)^{1/2} \leq \alpha_j(Q) \left(\frac{1}{\mu(Q)}\int_{Q} |f|^2  d\mu \right)^{1/2}$$
and
$$\forall j\geq 0 \quad \left(\frac{1}{\mu(2^{j+1}Q)} \int_{S_j(Q)} \left| f-B_Q(f)\right|^2 d\mu \right)^{1/2} \leq \alpha_j(Q) \left(\frac{1}{\mu(Q)}\int_{Q} |f|^{p_0} d\mu \right)^{1/p_0},$$
with coefficients $\alpha_j(Q)$ satisfying
$$\sup_{Q \textrm{ ball}} \sum_{j\geq 0} \frac{\mu(2^{j+1}Q)}{\mu(Q)} \alpha_j(Q) <\infty.$$
Let $\omega$ be a weight. Then for all exponents $p\in {\mathcal W}_{\omega}(p_0,2)$, there exists a constant $C$ such that
$$ \forall f\in L^2 \cap L^p(\omega), \qquad  \|T(f)\|_{p,\omega d\mu} \leq C \|f\|_{p,\omega d\mu}.$$
\end{thm}

\deme We have proved in Theorem \ref{theo4}, that under these assumptions the maximal operator $M_{(p_0)'}$ is bounded by the operator of Hardy-Littlewood $M_{HL,2}$. Then the desired result is a consequence of Theorem \ref{poidsthm} with $\sigma=(p_0)'$.

\begin{rem}Like the comparison between Theorem \ref{theo4} and Theorem \ref{theo5},
 we can compare Theorem \ref{theo4poids} with Theorem 8.8 of \cite{AM} : Theorem \ref{theo4poids}
 needs simpler and weaker assumptions as far as the continuity on $L^p$ with any $p\in {\mathcal W}_{\omega}(p_0,2)$ is concerned.
\end{rem}

\section{Embedding of $H^1_{\epsilon, mol}$ into $L^1$.}
\label{particular}

Here we discuss conditions on ${\mathbb B}=(B_Q)_Q$ insuring the embedding of our Hardy spaces
 into $L^1$. We assume throughout this section that ${\mathbb B}$ satisfies some decay estimates : for $M''$ a large enough exponent, there exists a constant $C$ such that \be{decay} \forall i\geq 0,\
\forall k\geq 0,\ \forall f\in L^2,\ \textrm{supp}(f) \subset 2^i Q \qquad \left\| B_Q(f) \right\|_{2,S_k(2^{i}Q)} \leq C 2^{-M''k}
\|f\|_{2,2^{i}Q}. \ee
In the sequel, all results about the atomic space $H^1_{ato}$ only require (\ref{decay}) with $i=0$. If we want to work with the molecular space $H^1_{\epsilon,mol}$, then we require $(\ref{decay})$ for all $i\geq 0$. We have the following imbedding~:

\begin{prop} \label{contL1} We have the following inclusions~:
$$ \forall \epsilon>0, \qquad H^{1}_{ato} \hookrightarrow H^{1}_{\epsilon, mol} \hookrightarrow L^1.$$
\end{prop}

\deme We claim that all $\epsilon$-molecules (and atoms) are bounded in $L^1$. In fact, using (\ref{decay})
\begin{align*}
\left\| B_Q(f_Q) \right\|_{1} & \leq  \sum_{i\geq 0} \left\| B_Q(f_Q{\bf 1}_{S_i(Q)}) \right\|_{1} \leq \sum_{i\geq 0} \sum_{k\geq 0} \left\| B_Q(f_Q{\bf 1}_{S_i(Q)}) \right\|_{1,S_k(2^{i}Q)} \\
 & \lesssim \sum_{i\geq 0} \sum_{k\geq 0} \mu(2^{i+k} Q)^{1/2} \left\| B_Q(f_Q{\bf 1}_{S_i(Q)}) \right\|_{2,S_k(2^{i}Q)} \\
 & \lesssim \sum_{i\geq 0} \sum_{k\geq 0} \mu(2^{i+k} Q)^{1/2} 2^{-kM''} \left\|f_Q \right\|_{2,S_i(Q)} \\
 & \lesssim \sum_{i\geq 0} \sum_{k\geq 0} \mu(2^{i+k} Q)^{1/2} 2^{-M''k} \mu(2^iQ)^{-1/2} 2^{-\epsilon i} \\
 & \lesssim \sum_{i\geq 0} \sum_{k\geq 0} 2^{k\delta/2} 2^{-M''k} 2^{-\epsilon i} \lesssim 1.
 \end{align*}
 Here we use the estimates for $f_Q$, the doubling property of $\mu$ and the fact that $M''$ is large enough ($M''>\delta/2$ works). So we obtain that all $\epsilon$-molecules are bounded in $L^1$, and we can deduce the imbedding using Lemma \ref{lemdense2}. \findem

\mb We recall that with this inclusion, in the study of the interpolation problem, we have obtained a more precise result (Corollary \ref{corint}).

\begin{cor} The spaces $H^1_{ato}$ and $H^1_{\epsilon, mol}$ are Banach spaces.
\end{cor}

\deme It is obvious that these spaces are normed vector spaces. We must verify the
 completeness. The proof is easy by
using the following well-known condition : for $\epsilon>0$,
$H^1_{\epsilon,mol}$ is a Banach space if for all sequences
$(h_i)_{i\in\N}$ of $H^1_{\epsilon,mol}$ satisfying
$$\sum_{i\geq 0} \left\| h_i \right\|_{H^1_{\epsilon,mol}} <\infty, $$
the series $\sum h_i$ converges in the Hardy space
$H^1_{\epsilon,mol}$. This is true because each molecular
decomposition is absolutely convergent in $L^1$-sense so we can
define the series $\sum_{i} h_i$ as a measurable function in $L^1$
and it is easy to prove the convergence for the $H^1_{\epsilon,mol}$
norm. \findem

\mb We have another corollary~:

\begin{cor} We also have the inclusions~:
$$L^\infty \subset (H^1_{\epsilon ,mol})^{*} \subset (H^1_{ato})^{*}.$$
\end{cor}

\mb  Now we want to compare our abstract Hardy spaces $H^1_{ato}$ and $H^1_{\epsilon, mol}$ with the classical Hardy space $H^1_{CW}$ of Coifman-Weiss. To do this we must define for a ball $Q$ the function $A_Q^{*}({\bf 1}_X)$.
Let $m:=B_Q(f_Q)$ an $\epsilon$-molecule of $H^1_{\epsilon, mol}$. We have seen that the integral
$$ \int_{X} A_Q(f_Q)(x) d\mu(x) $$
converges in $L^1$ sense (due to the decay of $f_Q$ and the "off-diagonal" decay of $A_Q$). Also we get
$$ \forall f\in L^1(X), \qquad \left|\int_{X} A_Q(f)(x) d\mu(x) \right| \lesssim  \sup_{i\geq 0} \|f\|_{2,S_i(Q)} \left(\mu(2^{i}Q)\right)^{1/2} 2^{\epsilon i}.$$
We can also consider the linear functional $f \rightarrow  \int_{X} A_Q(f)(x) d\mu(x)$, which will be denoted $A_Q^*({\bf 1}_{X})$, defined on the space
$$ Mol_{\epsilon,Q}:=\left\{ f\in L^1(X),\ \|f\|_{Mol_{\epsilon,Q}}:=\sup_{i\geq 0} \|f\|_{2,S_i(Q)} \left(\mu(2^{i}Q)\right)^{1/2} 2^{\epsilon i}<\infty\right\}.$$
With this preparation we have the following comparison between
Hardy spaces~:

\begin{prop} \label{inclus} Let $\epsilon\in]0,\infty]$. The inclusion $H^1_{\epsilon, mol} \subset H^1_{CW}$ is equivalent to the fact that for all $Q\in {\mathcal Q}$, $(A_Q)^* ({\bf 1}_X)={\bf 1}_{X}$ in $(Mol_{\epsilon,Q})^{*}$.
In this case for all $\epsilon'\geq \epsilon$ we have the inclusions $H^1_{ato} \subset H^1_{\epsilon',mol} \subset H^1_{\epsilon, mol} \subset H^1_{CW}$.
\end{prop}

\deme \\
$1-)$ Proof of the sufficiency. \\
Let $m:=B_Q(f_Q)$ an $\epsilon$-molecule of $H^1_{\epsilon, mol}$. First we want to prove that the integral of $m$ is equal to $0$. To show this we use the definition $B_Q=Id-A_Q$. By definition we have
$$ \int_{X} A_Q(f_Q)(x) d\mu(x)  = \langle A_Q(f_Q), {\bf 1}_{X} \rangle = \langle f_Q, A_Q^*({\bf 1}_{X}) \rangle.$$
So by the assumption, we have
$$\int_X A_Q(f_Q) d\mu = \int_X f_Q d\mu.$$
So we have shown that $\int B_Q(f_Q) d\mu= 0$. \\
Now we will prove that the molecule $m$ satisfies the good decay
around the ball $Q$ and so is an $H^1_{CW}$ molecule associated to
the ball $Q$. With the previous "off-diagonal" decay of $B_Q$, the
uniform $L^2$-boundedness of $B_Q$ and the assumptions over $f_Q$,
we have that for all $j\geq 0$ with $M''$ large enough~:
\begin{align*}
\left\| B_Q(f_Q) \right\|_{2,S_j(Q)} & \leq \sum_{i\geq 0} \left\| B_Q(f_Q{\bf 1}_{S_i(Q)}) \right\|_{2,S_j(Q)} \\
 & \leq \sum_{k=0}^{j} \left\| B_Q(f_Q{\bf 1}_{S_k(Q)}) \right\|_{2,S_j(Q)} + \sum_{k=j}^{\infty} \left\| B_Q(f_Q{\bf 1}_{S_k(Q)})
 \right\|_{2,X}\\
 & \leq \sum_{k=0}^{j} \left\| B_Q(f_Q{\bf 1}_{S_k(Q)}) \right\|_{2,S_{j-k}(2^{k}Q)} + \sum_{k=j}^{\infty} \left\|f_Q\right\|_{2,S_k(Q)} \\
 & \lesssim \sum_{k=0}^j 2^{-M''(j-k)} \|f_Q\|_{2,S_k(Q)} + \sum_{k=j}^{\infty} \mu(Q)^{-1/2}2^{-k \epsilon}\\
 & \lesssim \sum_{k=0}^j 2^{-M''(j-k)} \mu(Q)^{-1/2} 2^{-k\epsilon} + \mu(Q)^{-1/2}2^{-j\epsilon} \\
 & \lesssim \mu(Q)^{-1/2} 2^{-j\epsilon}.
\end{align*}
So we have proved that $m$ satisfies the decay of a classical
molecule in the space $H^1_{CW}$. So $m\in
H^1_{CW}$ and then $H^1_{\epsilon, mol} \subset H^1_{CW}$ follows from Lemma \ref{lemdense2}. \\
$2-)$ Proof of the necessity.\\
We assume that $H^1_{\epsilon, mol}$ is included in $H^1_{CW}$. Let $Q$ be a ball. So for each $\epsilon$-molecule $m=B_Q(f_Q)$ we have that
$$ \int B_Q(f_Q) d\mu = \int f_Q d\mu(Q) - \int A_Q(f_Q) d\mu =0.$$
We can compute the difference of the two integrals, because we have seen that they converge in $L^1$ sense. So by definition we have that for all $f\in Mol_{\epsilon,Q}$
$$ \int f_Q d\mu = \langle A_Q^*({\bf 1_X}), f \rangle,$$
which means that $A_Q^*({\bf 1}_X)={\bf 1}_{X}$ in $(Mol_{\epsilon,Q})^{*}$. \findem

\begin{rem} Under the assumption of the previous proposition and Theorem \ref{theogeneh}, it is interesting to note that the space $H^1_{ato}$ is smaller than the Hardy space $H^1_{CW}$ still it is big enough in $L^1$ to get the $L^p$ spaces by interpolation with $L^2$.
\end{rem}

\begin{rem} In \cite{DY} (see Subsection \ref{elliptic2}) the authors study whether their space $BMO_L$ satisfies $BMO \subset BMO_L$. In Proposition 6.7 of \cite{DY}, they prove that the inclusion $BMO \subset BMO_{L}$ is equivalent to the fact that for all $r>0$, $e^{-rL}({\bf 1}_{\R^n})={\bf 1}_{\R^n}$. By Theorem 3.1 of \cite{DY}, we know that $(H^1_{L})^{*} = BMO_{L^*}$. Also we have that
$$ H^1_{L} \subset H^1_{CW} \Longleftrightarrow  \ \forall r>0,\ e^{-rL^*}({\bf 1}_{\R^n})={\bf 1}_{\R^n}.$$
We recall (Subsection \ref{elliptic2}) that with the choice
$$ B_Q(f) := f-e^{-r_Q^2L}(f)$$
we have the inclusion
$$ \forall \epsilon>0, \qquad H^1_{\epsilon, mol} \subset H^1_{L}.$$
The previous proposition shows that we have~:
$$ H^1_{\epsilon, mol} \subset H^1_{CW} \Longleftrightarrow  \ \forall r>0,\ e^{-rL^*}({\bf 1}_{\R^n})={\bf 1}_{\R^n} \Longleftrightarrow H^1_{L} \subset H^1_{CW}.$$
\end{rem}

\section{An application to maximal $L^q$ regularity on Lebesgue spaces.}
\label{seven}

In this section, we apply the previous general theory to maximal $L^q$ regularity for Cauchy Problem. First we recall the subject of maximal $L^q$ regularity.

\mb Let $(Y,d_Y,\nu)$ a space of homogeneous type. Let $L$ the infinitesimal generator of an analytic semigroup of operators on $L^p(Y)$ and $J=(0,l],$ $l>0$ or $J=(0,+\infty)$ (in the second case, one has to assume that $A$ generates a bounded analytic semigroup).

\gb Consider the Cauchy problem
$$\left\{\begin{array}{ll}
\frac{du}{dt}(t)-Lu(t)=f(t), &  t\in J,\\
u(0)=0, &
\end{array}\right.
$$
where $f: J\rightarrow X$ is given.
If $e^{tL}$ is the semigroup generated by $L$, $u$ is given by
$$u(t)=\int_0^t e^{(t-s)L}f(s)ds.$$
For fixed $q \in (1,+\infty)$, one says that there is maximal $L^q$ regularity on $L^p(Y)$ for the problem if for every $f \in
L^q(J, L^p(Y))$, $\frac{\partial u}{\partial t}$(or $Lu$) belongs to $L^q(J, L^p(Y))$. For the maximal $L^q$
regularity, we refer to \cite{CV}, \cite{CD}, \cite{CD1}, \cite{DS}, \cite{HP}, \cite{L} etc.

\mb We now define an operator $T$~:
\begin{defn} With $L$ the generator, we define the operator~:
$$ Tf(t,x)=\int_0^t  \left[Le^{(t-s)L} f(s,.) \right](x) ds. $$
\end{defn}

\mb Let $p,q\in]1,\infty[$ be two exponents. We know that the maximal $L^q$ regularity on $L^p(Y)$ is equivalent to the fact that $T$ is bounded on $L^p(J \times Y)$. That is why we study this operator.
Of course, the problem of maximal $L^q$ regularity is completely understood by the abstract result in \cite{weis}. Here we want to remain as concrete as possible. In particular case we will see that the $H^1_{\epsilon,mol}-L^1$ continuity of the operator $T$ below depends only on $L^2$ assumptions. It is only when we want to deduce $L^p$ estimates that we need stronger assumptions which imply $R$-boundedness used in \cite{weis}.

\gb We need some conditions on our semigroup $e^{tL}$.

\begin{defn} Let ${\mathcal T}:=(T_t)_{t\in J}$ be a collection of $L^2(Y)$-bounded operators. We will say that ${\mathcal T}$ satisfies  off-diagonal $L^2-L^2$ estimates at ``low scale'' if there exists a function $\gamma$ satisfying
\be{gammaint} \forall p\geq 0, \qquad \sup_{u\geq 1} \gamma(u) u^{p} <\infty , \ee
such that for all balls $B\subset Y$ of radius $r$, for all functions $f$ supported in $B$ if $u\leq 2^{j+1}r$ then
\be{offdiagonallow}
 \left( \frac{1}{\nu(2^{j+1}B)} \int_{S_j(B)} \left|T_{u^2} (f) \right|^{2} d\nu \right)^{1/2} \leq \frac{\nu(B)}{\nu(2^{j+1}B)} \gamma\left(\frac{2^{j+1}r}{u} \right) \left( \frac{1}{\nu(B)}\int_{B} |f|^2 d\nu \right)^{1/2}.\ee
We also write ${\mathcal T} \in {\mathcal O}_{low}(L^2-L^2)$. \\
We will say that ${\mathcal T}$ satisfies off-diagonal $L^2-L^2$ estimates at ``high scale'' if for all balls $B\subset Y$ of radius $r$, for all functions $f$ supported in $B$
\be{offdiagonalhigh}  u\geq 2^{j+1}r \Longrightarrow
 \left( \frac{1}{\nu(2^{j+1}B)} \int_{S_j(B)} \left|T_{u^2} (f) \right|^{2} d\nu \right)^{1/2} \leq \frac{\nu(B)}{\nu(\frac{u}{r}B)} \left( \frac{1}{\nu(B)}\int_{B} |f|^2 d\nu \right)^{1/2}.\ee
We write also ${\mathcal T} \in {\mathcal O}_{high}(L^2-L^2)$. \\
In these two cases, the ``scale'' corresponds to the ratio between the parameter $u$ and the size $2^{j}r$ of the corona $S_j(B)$, where we estimate the operator $T_{u^2}$. \\
Let $q_0\in[2,\infty]$ be a fixed exponent, we say that ${\mathcal T}$ satisfies weak off-diagonal $L^2-L^{q_0}$ estimates if there exist coefficients $(\beta_j)_{j\geq 0}$ satisfying
\be{beta} \sum_{j\geq 0} 2^{j} \beta_j <\infty \ee
such that for all balls $B$ and for all functions $f\in L^2(Y)$ we have
 \be{offdiag} \left( \frac{1}{\nu(B)} \int_B \left|T_{r_Q^2} (f) \right|^{q_0} d\nu \right)^{1/q_0} \leq \sum_{j\geq 0} \beta_{j} \left( \frac{1}{\nu(2^{j} B)}\int_{2^{j}B} |f|^2 d\nu \right)^{1/2}.\ee
We also  write ${\mathcal T} \in {\mathcal O}_{weak}(L^2-L^{q_0})$.
\end{defn}

\begin{rem}

\mb $1-)$ It is easy to check that if a collection ${\mathcal T} \in {\mathcal O}_{high} (L^2-L^2)$ then ${\mathcal T}^*:=(T_t^*)_{t>0} \in {\mathcal O}_{weak}(L^2-L^2)$ with some coefficients $\beta_j$ satisfying
$$ \forall j\geq 0, \qquad \beta_j \leq \gamma(2^{j}).$$
$2-)$ All these conditions are satisfied with $q_0=\infty$ if the kernel $K_t$ of the operator $T_{t}$ admits some gaussian estimates like
$$ \left|K_t(x,y)\right| \lesssim \frac{1}{\nu(B(x,t^{1/2}))} e^{-d(x,y)^2/t}.$$
\end{rem}

\mb The main result of this section is the following one~:

\begin{thm} \label{thresume} Let $L$ be a generator of a bounded analytic semigroup ${\mathcal T}:=(e^{tL})_{t>0}$ on $L^2(Y)$ such that
${\mathcal T},(tLe^{tL})_{t>0}$ and $(t^2L^2e^{tL})_{t>0}$ belong to ${\mathcal O}_{low}(L^2-L^2) \cap {\mathcal O}_{high}(L^2-L^2)$. Then there exists a collection ${\mathbb B}=(B_Q)_{Q}$ of $L^2(J \times Y)$-bounded operators such that for all $\epsilon>0$ the operator $T$ is continuous from $H^1_{\epsilon, mol}(J \times Y)$ to $L^1(J\times Y)$. \\
In addition if ${\mathcal T}^*:=(e^{tL^*})_{t>0} \in {\mathcal O}_{weak}(L^2-L^{q_0})$ for a $q_0\in[2,\infty]$ then the maximal operator $M_{q_0}$ (defined by (\ref{opeM})) is bounded by the Hardy-Littlewood maximal operator $M_{HL,2}$.
\end{thm}

\mb We separate the proof in several steps. First we are going to
describe the choice of the collection ${\mathbb B}$. Then we will
check that the assumption (\ref{operh}) and the one about $M_{q_0}$
are satisfied. To finish the proof, we will show the $H^1_{\epsilon,
mol}-L^1$ boundedness of $T$ in Theorem \ref{theo2ap}.

\gb Equip $X=J \times Y$ with the parabolic quasi-distance $d$ and the measure $\mu$ defined by~:
$$d \Big( (t_1,y_1),(t_2,y_2) \Big) = \max\left\{ d_Y(y_1,y_2) , \sqrt{|t_1-t_2|} \right\} \qquad \textrm{and} \qquad  d\mu = dt \otimes d\nu.$$
If we write $\delta$ for the homogeneous dimension of the space $(Y,d_Y,\nu)$, then the space $X$ is of homogeneous type with homogeneous dimension $\delta+2$. We explain how to choose the collection $(B_Q)_{Q\in \mathcal{Q}}$ in this special case. We choose $\varphi \in \s(\R^+)$ such that $ \int_{ \R^+} \varphi(t)dt=1$ and $\varphi(t):=0$ for all $t<0$ ($\varphi$ does not need to be continuous at $0$). In fact we shall use only the fast decay of $\varphi$ and we will never consider regularity about it. In addition, we have added a condition for the support. This is a "physical" heuristids : this condition permits to define $A_Q(f)(t,x)$ by (\ref{opA}) with only $(f(\sigma,y))_{\sigma\leq t}$, which corresponds to the "past informations" about $f$. However we do not really need this assumption in the sequel. \\
For each cube $Q$ of $X$, we
write $r_Q$ its radius and we define the $B_Q$ operator as~:
\be{bq} B_Q=B_{r_Q^2} \qquad \textrm{with} \qquad B_r(f):= f-A_r(f), \ee
where the operator $A_r$ is defined by~: \be{opA} A_{r} (f)(t,x) :=\int _{\sigma =0}^{+\infty} \varphi_{r}(t-\sigma)
e^{rL}(f(\sigma,.))(x)  d\sigma. \ee Here we write $\varphi_{r}$ as the $L^1(\R)$ normalized function $\varphi_r(t):=r^{-1} \varphi(t/r)$. In fact, the integral for $\sigma\in[0,\infty[$ is reduced to $[0,t]$, due to the fact that $\varphi$ is supported in $\R^+$.

\mb Now we prove that these operators $B_Q$ satisfy the ``good''
conditions. First we have the assumption (\ref{operh})~:

\begin{prop} There is a constant $0<A'<\infty$ so that for all $r>0$ the operator $A_r$ is $L^2(X)$ bounded and we have~:
$$ \left\| A_r \right\|_{L^2 \rightarrow L^2} \leq A'.$$
\end{prop}

\deme By definition the semigroup $e^{rL}$ is $L^2(Y)$-bounded so we have the following estimates~:
\begin{align*}
\left\| A_{r} (f) \right\|_{2} & \leq \left\| \int _{\sigma =0}^{+\infty} \int _{y \in Y} \left|\varphi_{r}(t-\sigma)\right|
\left\|e^{rL}(f(\sigma,.)) \right\|_{2,d\nu}  d\sigma \right\|_{2,dt} \\
 & \lesssim \left\| \int _{\sigma =0}^{+\infty}  \left|\varphi_{r}(t-\sigma)\right|
\left\|f(\sigma,.) \right\|_{2,d\nu}  d\sigma \right\|_{2,dt} \\
 & \lesssim \left\| \int _{\sigma =0}^{+\infty}  \left|\varphi_{r}(\sigma)\right|
\left\|f(t-\sigma,.) \right\|_{2,d\nu}  d\sigma \right\|_{2,dt} \\
 & \lesssim \left\|\varphi_{r}\right\|_1 \left\|f \right\|_{2,d\mu} \lesssim \left\|f \right\|_{2,d\mu}.
\end{align*}
So we have proved that $A_r$ is $L^2(X)$-bounded and its boundedness is uniform for $r>0$.
\findem

\begin{thm} \label{assumcas} Define the maximal operator
$$ M_{q_0} (f)(\sigma,x):= \sup_{\genfrac{}{}{0pt}{}{Q \textrm{ ball}}{(\sigma,x)\in Q}} \left( \frac{1}{\mu(B)} \int_{B} \left| A_Q^*(f)\right|^{q_0} d\mu \right)^{1/q_0}.$$
If $(e^{tL^*})_{t>0} \in {\mathcal O}_{weak}(L^2-L^{q_0})$ then $M_{q_0}$ is bounded by the Hardy-Littlewood maximal operator $M_{HL,2}$ on $X$.
\end{thm}

\deme  Let $Q$ be a ball containing the point $(\sigma,x)\in X$ and $r_Q$ be its radius.
For $f,g\in L^2(X)$ we have~:
\begin{align*}
\langle A_{Q} (f) ,g \rangle & := \int_{(t,x)\in X} \int _{\sigma =0}^{+\infty} \varphi_{r_Q^2}(t-\sigma)
e^{r_Q^2 L}(f(\sigma,.))(x) g(t,x)  d\sigma dt d\nu(x) \\
 & = \int_{(t,x)\in X} \int _{\sigma =0}^{+\infty} \varphi_{r_Q^2}(t-\sigma)
f(\sigma,x) \left[\left(e^{r_Q^2 L}\right)^*g(t,.)\right](x)  d\sigma dt d\nu(x).
\end{align*}
So we conclude that~: \be{Aetoile} A_Q^*(g)(\sigma,x):= \int_{t\in
\R^+}  \varphi_{r_Q^2}(t-\sigma) \left[\left(e^{r_Q^2
L}\right)^*g(t,.)\right](x) dt.\ee By using the Minkowski
inequality, we also have  that
$$ \left ( \int_{Q} \left| A_Q^*(g)\right|^{q_0} d\mu \right)^{\frac{1}{q_0}}  \leq \int_{t\in \R^+}  \left\| \varphi_{r_Q^2}(t-\sigma) \left[\left(e^{r_Q^2 L}\right)^*g(t,.)\right](x) \right\|_{q_0,d\nu(x)d\sigma} dt. $$
By definition of the parabolic quasi-distance, we can write
$$ Q = I \times B$$
with $I$ an interval of lenght $r_Q^2$ and $B$ a ball of $Y$ of radius $r_Q$. Then we have~:
\begin{align*}
\lefteqn{ \left( \int_{Q} \left| A_Q^*(f)\right|^{q_0} d\mu \right)^{1/q_0}  \leq} & & \\
& &  \int_{t\in \R^+}  \left\| \varphi_{r_Q^2}(t-\sigma){\bf 1}_{I}(\sigma) \right\|_{q_0,d\sigma} \left\|{\bf 1}_{B}(x) \left(e^{r_Q^2 L}\right)^*g(t,.)(x) \right\|_{q_0,d\nu(x)} dt.
\end{align*}
With the assumption (\ref{offdiag}), we obtain
\begin{align*}
\lefteqn{ \left( \int_{Q} \left| A_Q^*(f)\right|^{q_0} d\mu \right)^{1/q_0} \leq } & & \\
 & &   \sum_{j\geq 0} \int_{t\in \R^+}  \left\| \varphi_{r_Q^2}(t-\sigma) {\bf 1}_{I}(\sigma)\right\|_{q_0,d\sigma}  \beta_j \frac{\nu(B)^{1/q_0}}{\nu(2^{j}B)^{1/2}} \left\|g(t,x) {\bf 1}_{2^{j}B}(x)\right\|_{2,d\nu(x)} dt.
 \end{align*}
Now we decompose the integration over $t$ by~:
\begin{align*}
\lefteqn{ \left( \int_{Q} \left| A_Q^*(f)\right|^{q_0} d\mu \right)^{1/q_0}  \leq} & & \\
 & & \sum_{j\geq 0} \sum_{k\geq 0} \int_{t \in S_k(I)}  \left\| \varphi_{r_Q^2}(t-\sigma) {\bf 1}_{I}(\sigma) \right\|_{q_0,d\sigma}  \beta_j \frac{\nu(B)^{1/q_0}}{\nu(2^{j}B)^{1/2}} \left\|g(t,x) {\bf 1}_{2^jB}(x) \right\|_{2,d\nu(x)} dt.
 \end{align*}
With the Cauchy-Schwarz inequality, we have
\begin{align*}
\lefteqn{\left( \int_{Q} \left| A_Q^*(f)\right|^{q_0} d\mu \right)^{1/q_0} } & & \\
& & \lesssim \sum_{j\geq 0} \sum_{k\geq 0} r_Q^{-2}
\left(1+ 2^{k}\right)^{-l} r_Q^{2/q_0} \beta_j \frac{\nu(B)^{1/q_0}}{\nu(2^{j}B)^{1/2}} (2^{k} r_Q^2)^{1/2}\left\|g(t,x){\bf 1}_{2^kI \times 2^jB}(t,x) \right\|_{2,dt d\nu(x)} \\
& & \lesssim \sum_{j\geq 0} \sum_{k\geq 0} r_Q^{-1+2/q_0}
\left(1+ 2^{k}\right)^{-l+1/2} \beta_j \frac{\nu(B)^{1/q_0}}{\nu(2^{j}B)^{1/2}} \left\|g(t,x) {\bf 1}_{2^kI \times 2^jB}(t,x)\right\|_{2,dtd\nu(x)}.
\end{align*}
Here $l$ is an integer as large as we want, due to the fast decay of $\varphi$. Using the maximal Hardy-Littlewood operator, we have
$$ \left\|g(t,x) {\bf 1}_{2^kI \times 2^jB}(t,x) \right\|_{2,dtd\nu(x)} \leq \mu\left(\max\{2^j,2^{k/2}\} Q\right)^{1/2}  \inf_{Q} M_{HL,2}(g).$$
So we obtain
\begin{align*}
\lefteqn{\left( \int_{Q} \left| A_Q^*(g)\right|^{q_0} d\mu \right)^{1/q_0} \leq} & & \\
& &  \left[ \sum_{j\geq 0} \sum_{k\geq 0} r_Q^{-1+2/q_0} \left(1+ 2^{k}\right)^{-l+1/2} \beta_j \frac{\nu(B)^{1/q_0}}{\nu(2^{j}B)^{1/2}}  \mu\left(\max\{2^j,2^{k/2}\} Q\right)^{1/2} \right]  \inf_{Q} M_{HL,2}(g). \label{equ1}
\end{align*}
We now estimate the sum over the parameters $j$ and $k$. We have the two following cases.
Write
$$  S_1 := \sum_{j\geq k/2 \geq 0} r_Q^{-1+2/q_0} \left(1+ 2^{k}\right)^{-l+1/2} \beta_j \frac{\nu(B)^{1/q_0}}{\nu(2^{j}B)^{1/2}}  \mu\left( 2^{j} Q\right)^{1/2} $$
and
$$ S_2 := \sum_{k/2 \geq j \geq 0} r_Q^{-1+2/q_0} \left(1+ 2^{k}\right)^{-l+1/2} \beta_j \frac{\nu(B)^{1/q_0}}{\nu(2^{j}B)^{1/2}}  \mu\left(2^{k/2} Q\right)^{1/2} .$$
We must estimate these two sums. For the first, we use that $\mu(Q)= |I| \nu(B)=r_Q^2 \nu(B)$ to have
\begin{align*}
  S_1 & \leq  \sum_{j\geq k/2 \geq 0} 2^{j} \left(1+ 2^{k}\right)^{-l+1/2} \beta_j \frac{\mu(Q)^{1/q_0}}{\mu(2^{j}Q)^{1/2}}  \mu\left( 2^{j} Q\right)^{1/2} \\
   & \leq  \mu(Q)^{1/q_0} \sum_{j\geq k/2 \geq 0} 2^{j} \left(1+ 2^{k}\right)^{-l+1/2} \beta_j \\
   & \leq  \mu(Q)^{1/q_0} \sum_{j\geq 0} 2^{j} \beta_j \lesssim \mu(Q)^{1/q_0}. \\
\end{align*}
In the last inequality, we have used the assumption (\ref{beta}) about the coefficients $(\beta_j)_j$. \\
For the second sum, we have (with the doubling property of $\mu$ and $l$ large enough)
 \begin{align*}
  S_2 & \leq r_Q^{2/q_0} \nu(B)^{1/q_0} \sum_{k/2 \geq j \geq 0} r_Q^{-1} \left(1+ 2^{k}\right)^{-l+1/2} \beta_j \left(\frac{ \mu(2^{k/2} Q) }{\nu(2^{j}B)}\right)^{1/2}  \\
   & \lesssim \mu(Q)^{1/q_0} \sum_{k/2 \geq j \geq 0}  r_Q^{-1} \left(1+ 2^{k}\right)^{-l+1/2} \beta_j \left(\frac{ \mu(2^{j} Q) }{\nu(2^{j}B)}\right)^{1/2} 2^{(k/2-j)(\delta+2)/2} \\
   & \lesssim \mu(Q)^{1/q_0} \sum_{k/2 \geq j \geq 0}  \left(1+ 2^{k}\right)^{-l+1/2} \beta_j 2^{j} 2^{(k/2-j)(\delta+2)/2}  \\
   & \lesssim \mu(Q)^{1/q_0} \sum_{j \geq 0}  \left(1+ 2^{j}\right)^{-l+4+\delta/2} \beta_j 2^{-j(\delta/2+1)} \lesssim \mu(Q)^{1/q_0}.
\end{align*}
So we have proved that there exists a constant $C$ (independant on $g$ and $Q$) such that~:
$$ \left( \int_{Q} \left| A_Q^*(g)\right|^{q_0} d\mu \right)^{1/q_0} \leq C \mu(Q)^{1/q_0}  \inf_{Q} M_{HL,2}(g). $$
 We can also conclude that
 $$M_{q_0}(f) \lesssim M_{HL,2}(g).$$
  \findem

\mb The assumption on $M_{q_0}$ satisfied : $M_{q_0}$ is bounded by
$M_{HL,2}$. To apply the previous abstract result about
interpolation, we must now show that our operator $T$ is continuous
from our Hardy spaces into $L^1$ and prove the $L^2$-boundedness of
$T$.

\begin{thm} $T$ is $L^2(X)$-bounded.
\end{thm}

\mb This fact was proved in \cite{DS} because it is equivalent to the maximal $L^2$ regularity on $L^2(Y)$.

\mb Now we want to apply the general Theorem \ref{theo2h}. We have the following result~:

\begin{thm} \label{theo2ap}  Let $L$ be a generator of a bounded analytic semigroup on $L^2(Y)$. Assume that $(e^{tL})_{t>0},(tLe^{tL})_{t>0}$ and $(t^2L^2e^{tL})_{t>0}$ belong to ${\mathcal O}_{low}(L^2-L^2) \cap {\mathcal O}_{high}(L^2-L^2)$. Then there exist coefficients $\alpha_{j,k}$ such that for all balls $Q\subset X$, for all $k\geq 0,j\geq 2$ and for all functions $f$ supported in $S_k(Q)$
\be{theo2a} \left(\frac{1}{\mu(2^{j+k+1}Q)} \int_{S_j(2^kQ)} \left| T(B_Q(f))\right|^2 d\mu \right)^{1/2} \leq \alpha_{j,k} \left(\frac{1}{\mu(2^{k+1}Q)}\int_{S_k(Q)} |f|^2 d\mu \right)^{1/2}. \ee
In addition the coefficients $\alpha_{j,k}$ (independent in $Q$) satisfy \be{hypap} \Lambda:= \sup_{Q} \ \sup_{k\geq 0}
\left[\sum_{j\geq 2} \frac{\mu(2^{j+k+1}Q)}{\mu(2^{k+1}Q)} \alpha_{j,k} \right] <\infty. \ee
With Theorem \ref{theo2h}, these estimates imply the $H^1_{\epsilon, mol}(X)-L^1(X)$ boundedness of $T$ for all $\epsilon>0$.
\end{thm}

\deme We write $r=r_Q$ and $(t_0,x_0)$ the radius and the center of the ball $Q$ so we have defined $B_Q$ as $B_{r^2}$. The function $f$ is fixed. The parameter $j$ and $k$ are fixed too. We write $Q$ as the product $Q=I \times B$ with $I$ an interval of length $r_Q^2$ and $B$ a ball of $Y$ of radius $r_Q$.
We have
\begin{align*}
T B_{r^2}(f) (t,x) & = T(f)(t,x) - T A_{r^2}(f) (t,x ) \\
 & =\int _{0} ^{t} \left[Le^{(t-s)L}f(s,.)\right](x) ds-\int_{0}^{t} \left[L e^{(t-s)L} A_{r^2} f(s,.)\right](x) ds,
\end{align*}
where
\begin{align*}
 \left[L e^{(t-s)L} A_{r^2}f(s,.) \right](x) = L e^{(t-s)L} \left[\int_{\sigma=0}^{+\infty} \varphi
_{r^2} (s-\sigma)  e^{r^2 L} f(\sigma,.)  d\sigma \right](x).
\end{align*}
So we obtain
\begin{align}
\lefteqn{T(B_{r^2}f)(t,x)  =} & & \nonumber \\
 & & \int_\R \int_\R \varphi_{r^2}(s-\sigma)\left[{\bf 1}_{0<\sigma\leq t} Le^{(t-\sigma)L}f(\sigma,.)(x) -{\bf 1}_{0<s\leq t} L e^{(t-s+r^2)L}f(\sigma,.)(x) \right] d\sigma ds . \label{importante}
 \end{align}
We have three time-parameters $\sigma,t$ and $s$. As in the case of Calder\'on-Zygmund operators, the difference within the two brackets is very important. This will allow us to obtain the necessary decay for the coefficients $\alpha_{j,k}$.
We decompose into two domains~:
$$ D_1:=\left\{ (\sigma,t,s),\ 0\leq \sigma \leq t\leq s \right\} \quad \textrm{and} \quad D_2:=\left\{ (\sigma,t,s),\ 0\leq s,\sigma \leq t \right\}.$$
For $i\in\{1,2\}$ we set $D_i(t):=\left\{(\sigma,s);\ (\sigma,t,s)\in D_i\right\}$ and
\begin{align*}
\lefteqn{ U_i(f)(t,x) :=} & & \\
 & &  \iint_{D_i(t)} \varphi_{r^2}(s-\sigma)\left[{\bf 1}_{0<\sigma\leq t} Le^{(t-\sigma)L}f(\sigma,.)(x) -{\bf 1}_{0<s\leq t} L e^{(t-s+r^2)L}f(\sigma,.)(x) \right] d\sigma ds.
 \end{align*}
As $\varphi$ is supported in $\R^+$, we have decomposed 
\be{decompositionTB} T(B_{r^2}f)(t,x) = \sum_{i=1}^{2} U_i(f)(t,x). \ee
If we do not want to use the condition of the support of $\varphi$, there is a third term which is estimated as the first one. \\
We begin the study when one of the two terms, in the square brackets, vanishes. The radius $r$ is fixed for all the proof and we set
$$ \chi_N(y):=\frac{1}{r^2}\left(1+ \frac{|y|}{r^2} \right)^{-N}.$$
 \\
$1-)$ First case : $(\sigma,s)\in D_1(t)$. \\
Here we have the following expression~:
\begin{align*}
U_1(f)(t,x) = -\int_{t}^\infty \int_{0}^{t} \varphi_{r^2}(s-\sigma) L e^{(t-\sigma)L} f(\sigma,.)(x) d\sigma ds.
\end{align*}
There is no "cancellation" so we can directly estimate it by using the fast decay of $\varphi$.
For $N$ a large enough integer
\begin{align*}
\left| U_1(f)(t,x) \right| & \lesssim  r^2\int_{0}^{t} \chi_N(t-\sigma) \left| L e^{(t-\sigma)L} f(\sigma,.)(x)\right| d\sigma.
\end{align*}
By definition of the parabolic quasi-distance,
$$(t,x) \in S_j (2^k Q) \Longleftrightarrow  \left\{ \begin{array}{l}
d_Y(x,x_0) \simeq 2^{k+j}r \\
|t-t_0| \leq ( 2^{k+j+1}  r)^2 \end{array} \right.
\textrm{ or }
\left\{\begin{array}{l}
d_Y(x,x_0) \lesssim 2^{k+j}r \\
|t-t_0| \simeq (2^{k+j+1} r)^2 \end{array} \right. .
$$
So, as $f$ is supported in $2^kQ$, we have
\be{III} \left(\frac{1}{\mu(2^{j+k+1}Q) }\int_{S_j(2^kQ)} \left| U_1(f)(t,x) \right|^2 d\nu(x)dt \right)^{1/2} \lesssim I+ II \ee
with
$$ I := \left( \frac{r^4}{\mu(2^{j+k+1}Q)} \int_{2^{2(j+k)}I} \left( \int_{2^{2k}I} \chi_N(t-\sigma) \left\| L e^{(t-\sigma)L} f(\sigma,.)\right\|_{2,S_j(2^kB)} d\sigma\right)^2 dt \right)^{1/2}$$
and
$$ II :=  \left(\frac{r^4}{\mu(2^{j+k+1}Q)} \int_{S_{2j}(2^{2k}I)} \left(\int_{2^{2k}I} \frac{ \left\| L e^{(t-\sigma)L} f(\sigma,.)\right\|_{2,2^{k+j+1}B}}{2^{2N(k+j)}} d\sigma\right)^2 dt \right)^{1/2}.$$

\noindent $*$ Study of $I$. \\
By using off-diagonal estimates $L^2-L^2$ (\ref{offdiagonallow}), we know that
\begin{align*}
 \lefteqn{\frac{1}{\nu(2^{j+k+1}B)^{1/2}} \left\| L e^{(t-\sigma)L} f(\sigma,.)\right\|_{2,S_j(2^kB)} \leq} & & \\
 & &  \frac{\nu(2^{k+1}B)}{\nu(2^{k+j+1}B)|t-\sigma|} \gamma\left(\frac{2^{j+k}r}{\sqrt{|t-\sigma|}}\right)  \left( \frac{1}{\nu(2^{k+1}B)} \int_{2^{k+1} B} |f(\sigma,.)|^2 d\nu\right)^{1/2}.
 \end{align*}
That is why, by using Cauchy-Schwarz inequality and the equality
$$ \mu(2^{j+k+1}Q) = \nu(2^{j+k+1}B) 2^{2(j+k)}r^2,$$
we estimate $I$ by the product
\begin{align*}
  \left( \frac{1}{2^{2(k+j)}} \int_{2^{2k+2j}I} \int_{2^{2k}I} \chi_{2N}(t-\sigma) \left(\frac{\nu(2^{k+1}B)}{\nu(2^{k+j+1}B)|t-\sigma|}\right)^2 \gamma\left(\frac{2^{j+k}r}{\sqrt{|t-\sigma|}}\right)^2  d\sigma dt \right)^{1/2}      \\
    \qquad  2^{k}r\left( \frac{1}{\mu(2^{k+1}Q)} \int_{2^{k+1} Q} |f|^2 d\mu\right)^{1/2}.
\end{align*}
Then we get
\begin{align*}
  I & \lesssim \frac{\nu(2^{k+1}B)}{\nu(2^{k+j+1}B)} \Bigg[ \frac{1}{|2^{2k+2j} I|} \int_{2^{2k+2j}I} \int_{2^{2k}I} \chi_{2N}(t-\sigma)  \frac{2^{2k}r^4}{|t-\sigma|^2}  \\
    & \hspace{2cm}  \gamma\left(\frac{2^{j+k}r}{\sqrt{|t-\sigma|}}\right)^2 d\sigma dt \Bigg]^{1/2} \left( \frac{1}{\mu(2^{k+1}Q)} \int_{2^{k+1} Q} |f|^2 d\mu\right)^{1/2} \\
  & \lesssim \frac{2^k \nu(2^{k+1}B)}{\nu(2^{k+j+1}B)} \left(\int_{0}^{2^{2(j+k)}} \left(1+v \right)^{-2N}  \frac{1}{v^2} \gamma\left(\frac{2^{j+k}}{\sqrt{v}}\right)^2 dv \right)^{1/2} \\
  &  \hspace{5cm} \left( \frac{1}{\mu(2^{k+1}Q)} \int_{2^{k+1} Q} |f|^2 d\mu\right)^{1/2} \\
  & \lesssim \frac{2^{k}\nu(2^{k+1}B)}{\nu(2^{k+j+1}B)} 2^{-j-k}\left(\int_{1}^{\infty} \left(1+2^{k+j}v^{-2} \right)^{-2N} \gamma(v)^2 v dv\right)^{1/2} \\
  &  \hspace{5cm} \left( \frac{1}{\mu(2^{k+1}Q)} \int_{2^{k+1} Q} |f|^2 d\mu\right)^{1/2}.
\end{align*}
$*$ Study of $II$. \\
In this case, we have $t\in S_{2j}(2^{2k}I)$ and $\sigma \in 2^{2k}I$, so
 $$ |t-\sigma|\simeq 2^{2(j+k)}r^2.$$
By using off-diagonal estimates (\ref{offdiagonalhigh}), we know that
\begin{align*}
 \lefteqn{\frac{1}{\nu(2^{j+k+1}B)^{1/2}} \left\| L e^{(t-\sigma)L} f(\sigma,.)\right\|_{2,2^{k+j+1}B} \lesssim \hspace{2cm} } & & \\
 & &  \hspace{2cm} \frac{1}{2^{2(j+k)} r^2}  \left( \frac{1}{\nu(2^{k+1}B)} \int_{2^{k+1} B} |f(\sigma,.)|^2 d\nu\right)^{1/2}.
 \end{align*}
So we obtain that
\begin{align*}
  II & \lesssim   \left(\frac{1}{2^{2(j+k)}r^2} \int_{t\in S_{2j}(2^{2k}I)} \int_{2^{2k}I} 2^{-4N(k+j)} \frac{2^{2k}r^2 }{2^{4(j+k)} r^4}  d\sigma dt \right)^{1/2} & \\
  & & \hspace{-3cm} \left( \frac{1}{\mu(2^{k+1}Q)} \int_{2^{k+1} Q} |f|^2 d\mu\right)^{1/2} \\
  & \lesssim  2^{-4j}2^{-2N{k+j}} \left( \frac{1}{\mu(2^{k+1}Q)} \int_{2^{k+1} Q} |f|^2 d\mu\right)^{1/2}.
\end{align*}

\gb We have also the following estimate
\begin{align*} \lefteqn{I+ II \lesssim \left( \frac{1}{\mu(2^{k+1}Q)} \int_{2^{k+1} Q} |f|^2 d\mu\right)^{1/2} } & & \\
 & &    \left[2^{-4j} 2^{-N(k+j)}  + \frac{\nu(2^{k+1}B)}{\nu(2^{k+j+1}B)} 2^{-j-k} \left(\int_{1}^\infty \left(1+2^{k+j}v^{-2} \right)^{-2N} \gamma(v)^2 v dv \right)^{1/2}\right].
 \end{align*}
With (\ref{III}), here we can choose
$$ \alpha_{j,k}= \left[2^{-N(k+j)}  +  \frac{\nu(2^{k+1}B)}{\nu(2^{k+j+1}B)}2^{-j-k}\left(\int_{1}^\infty \left(1+2^{j}v^{-2} \right)^{-N} \gamma(v)^2 vdv \right)^{1/2} \right]$$
for $N$ a large enough integer. \\
$2-)$ Last case for $(\sigma,s)\in D_2(t)$ : $0\leq s,\sigma \leq t$.\\
The relation (\ref{importante}) gives us that~:
$$ U_2(f)(t,x)  =
\int_0^t \int_0^t \varphi_{r^2}(s-\sigma)\left[Le^{(t-\sigma)L}f(\sigma,.)(x) - L e^{(t-s+r^2)L}f(\sigma,.)(x) \right] d\sigma ds. $$
Here we use the time regularity. We have~:
\begin{align*}
\left|Le^{(t-\sigma)L}f(\sigma,.)(x) - L e^{(t-s+r^2)L}f(\sigma,.)(x) \right| & = \left| \int_{t-s+r^2}^{t-\sigma} \frac{\partial L e^{zL}f(\sigma,.)(x)}{\partial z} dz \right| \\
 & = \left| \int_{t-s+r^2}^{t-\sigma} L^2 e^{zL}f(\sigma,.)(x) dz \right|.
\end{align*}
Then we repeat the same arguments as before~: \\
\be{IIII} \left(\frac{1}{\mu(2^{j+k+1}Q)}\int_{S_j(2^kQ)} \left| U_2(f)(t,x) \right|^2 d\nu(x)dt \right)^{1/2} \lesssim I+ II \ee
with
\begin{align*}
 \lefteqn{I := \Bigg( \frac{1}{\mu(2^{j+k+1}Q)} \int_{2^{2k+2j}I} \left( \int_{2^{2k}I} \int_0^t \chi_N(s-\sigma) \right.  \hspace{3cm} } & & \\
 & & \left. \left. \hspace{3cm} \int_{t-s+r^2}^{t-\sigma} \left\| L^2 e^{zL} f(\sigma,.)\right\|_{2,S_j(2^kB)} dz ds d\sigma\right)^2 dt \right)^{1/2}
 \end{align*}
and
\begin{align*}
\lefteqn{ II :=  \left( \frac{1}{\mu(2^{j+k+1}Q)} \int_{t\in S_{2j}(2^{2k}I)} \left(\int_{2^{2k}I} \int_0^t \chi_N(s-\sigma)  \right. \right. \hspace{3cm} } & & \\
 & & \left. \left. \hspace{3cm} \int_{t-s+r^2}^{t-\sigma} \left\| L^2 e^{zL} f(\sigma,.)\right\|_{2,2^{k+j}B} dz ds d\sigma\right)^2 dt \right)^{1/2}.
 \end{align*}
$*$ Study of $I$. \\
By using off-diagonal estimates (\ref{offdiagonallow}), we know that
\begin{align*}
 \lefteqn{\frac{1}{\nu(2^{j+k+1}B)^{1/2}}\left\| L^2 e^{zL} f(\sigma,.)\right\|_{2,S_j(2^kB)} \leq \hspace{2cm} } & & \\
 & &  \hspace{2cm} \frac{\nu(2^kB)}{\nu(2^{k+j+1}B)z^2} \gamma\left(\frac{2^{j+k}r}{\sqrt{z}} \right)  \left( \frac{1}{\nu(2^{k}B)} \int_{2^k B} |f(\sigma,.)|^2 d\nu\right)^{1/2}.
\end{align*}
So we obtain
\begin{align*}
 \lefteqn{ I \lesssim  \frac{\nu(2^kB)}{\nu(2^{k+j+1}B)} \left( \frac{1}{\mu(2^{k}Q)} \int_{2^k Q} |f|^2 d\mu\right)^{1/2} {\Bigg [} \frac{1}{2^{2(k+j)} r^2} \int_{2^{2k+2j}I}\int_{2^{2k}I}  } & & \\
  & &  \left( \int_{0}^{t}  \frac{1}{r^2} \left(1+\frac{|s-\sigma|}{r^2}\right)^{-N} \left|\int_{t-s+r^2}^{t-\sigma}  \frac{2^{k}r}{z^2} \gamma\left(\frac{2^{j+k}r}{\sqrt{z}}\right) dz\right| ds\right)^2 d\sigma dt {\Bigg ]} ^{1/2} .
 \end{align*}
We use the inequality
\begin{align}
 \left| \int_{t-s+r^2}^{t-\sigma}  \frac{1}{z^2} \gamma\left(\frac{2^{j+k}r}{\sqrt{z}}\right) dz \right| & \leq \| \gamma\|_\infty \left|\frac{1}{t-s+r^2} - \frac{1}{t-\sigma} \right| = \|\gamma \|_\infty \frac{|s-\sigma+r^2|}{|t-\sigma| |r^2 + t-s|}  \nonumber \\
 & \lesssim \|\gamma\|_\infty \frac{1+ \frac{|s-\sigma|}{r^2}}{2^{2(j+k)} r^2 \left(1 + \frac{|t-s|}{r^2} \right)} \lesssim \|\gamma\|_\infty \frac{\left(1+ \frac{|s-\sigma|}{r^2}\right)^2}{2^{4(j+k)} r^2}, \label{difference}
 \end{align}
to have
\begin{align*}
 I & \lesssim \frac{\|\gamma\|_\infty \nu(2^kB)2^{k}}{\nu(2^{k+j+1}B)2^{4(k+j)}r}  \left( \frac{1}{|2^{2k+2j} I|} \int_{2^{2k+2j}I}  2^{2k}r^2 dt \right)^{1/2} \left( \frac{1}{\mu(2^{k}Q)} \int_{2^k Q} |f|^2 d\mu\right)^{1/2} \\
  & \lesssim \frac{\nu(2^kB)}{\nu(2^{k+j+1}B)} 2^{-4j} \|\gamma\|_\infty \left( \frac{1}{\mu(2^{k}Q)} \int_{2^k Q} |f|^2 d\mu\right)^{1/2}.
\end{align*}
Here we can choose
$$ \alpha_{j,k}= \frac{\nu(2^kB)}{\nu(2^{k+j+1}B)} 2^{-4j} \|\gamma\|_\infty \simeq \frac{\mu(2^kQ)}{\mu(2^{k+j+1}Q)} 2^{-2j} \|\gamma\|_\infty .$$
$*$ Study of $II$. \\
In this case, we have $|t-\sigma| \simeq 2^{2(k+j)}r^2$.
By using off-diagonal estimates (\ref{offdiagonalhigh}), we know that for $z\geq r^2$
$$ \frac{1}{\nu(2^{j+k+1}Q)^{1/2}}\left\| L^2 e^{zL} f(\sigma,.)\right\|_{2,2^{k+j}B} \leq \frac{1}{z^2}
  \frac{\nu(2^{k}B)}{\nu(\frac{z^{1/2}}{r}B)} \left( \frac{1}{\nu(2^{k}B)} \int_{2^k B} |f(\sigma,.)|^2 d\nu\right)^{1/2}. $$
So we obtain
\begin{align*}
 II & \lesssim  2^{-j} \left( \int_{2^{2k+2j}I} \int_{2^{2k}I} \left( \int_{0}^{t} \chi_N(s-\sigma) \int_{t-s+r^2}^{t-\sigma}  \frac{1}{z^2} \frac{\nu(2^{k}B)}{\nu(\frac{z^{1/2}}{r}B)} dz ds\right)^2 d\sigma dt \right)^{1/2}  \\
  & \hspace{7cm}  \left( \frac{1}{\mu(2^{k}Q)} \int_{2^k Q} |f|^2 d\mu\right)^{1/2}.
 \end{align*}
We use the inequality (based on the doubling property of $\nu$) that for all $z\in [t-\sigma,t-s+r^2]$ we have
\begin{align*}
 \frac{1}{\nu(\frac{z^{1/2}}{r}B)} & \lesssim  \frac{1}{\nu(\frac{|t-\sigma|^{1/2}}{r} B)} \left(1+ \frac{|t-\sigma|^{1/2}}{z^{1/2}} \right)^\delta \lesssim \frac{1}{\nu(\frac{|t-\sigma|^{1/2}}{r} B)}  \left(1+ \frac{|t-\sigma|^{1/2}}{(t-s+r^2)^{1/2}} \right)^\delta \\
  & \lesssim \frac{1}{\nu(\frac{|t-\sigma|^{1/2}}{r} B)} \left(1+ \frac{1+\frac{|t-\sigma|^{1/2}}{r}}{1+\frac{|t-s|^{1/2}}{r}} \right)^\delta \\
  & \lesssim \frac{1}{\nu(\frac{|t-\sigma|^{1/2}}{r} B)} \left(1+\frac{|s-\sigma|^{1/2}}{r} \right)^\delta.
\end{align*}
So as $N$ is an integer as large as we want, we can estimate
\begin{align*} II & \lesssim  \frac{\nu(2^{k}B)}{2^{j} \nu(2^{j+k+1}B)} \left(\int_{2^{2k+2j}I} \int_{2^{2k}I} \left(\int_{0}^{t}  \chi_N(s-\sigma) \int_{t-s+r^2}^{t-\sigma}  \frac{1}{z^2} dz ds \right)^2 d\sigma dt \right)^{1/2}  \\
  & \hspace{7cm}  \left( \frac{1}{\mu(2^{k}Q)} \int_{2^k Q} |f|^2 d\mu\right)^{1/2}.
\end{align*}
Then we use (\ref{difference}) (with $\gamma$ equal to the constant function) to finally obtain (with an other exponent $N$)
\begin{align*} II & \lesssim  \frac{\nu(2^{k}B)}{\nu( 2^{j+k+1}B)2^{j} 2^{4(j+k)} r^2} \left( \int_{2^{2k+2j}I} \int_{2^{2k}I} \left( \int_{0}^{t} \chi_N(s-\sigma) ds \right)^2 d\sigma dt \right)^{1/2}  \\
  & \hspace{7cm}  \left( \frac{1}{\mu(2^{k}Q)} \int_{2^k Q} |f|^2 d\mu\right)^{1/2} \\
  & \lesssim \frac{\nu(2^{k}B)2^{2k} r^2}{\nu( 2^{j+k+1}B)2^{4(j+k)} r^2}\left( \frac{1}{\mu(2^{k}Q)} \int_{2^k Q} |f|^2 d\mu\right)^{1/2} \\
  & \lesssim \frac{\mu(2^{k}Q)}{\mu( 2^{j+k+1}Q)2^{2j+2k}}\left( \frac{1}{\mu(2^{k}Q)} \int_{2^k Q} |f|^2 d\mu\right)^{1/2}.
\end{align*}
So here we can choose
$$ \alpha_{j,k}= \frac{\mu(2^{k}Q)}{\mu( 2^{j+k+1}Q)2^{2j+2k}}. $$

\noindent $3-)$ End of the proof.\\
With the decomposition (\ref{decompositionTB}), we have proved in the two previous points that we have the estimate
(\ref{theo2a}) with the coefficients $\alpha_{j,k}$ satisfying
\begin{align*}
 \alpha_{j,k} \lesssim & \  2^{-N(k+j)}  + \frac{\nu(2^{k}B)}{\nu(2^{k+j+1}B)}2^{-j} \left(\int_{1}^\infty \left(1+2^{j}v^{-2} \right)^{-2N}\gamma(v)^2 vdv \right)^{1/2}  \\
 &  \ + \frac{\mu(2^{k}Q)}{\mu( 2^{j+k+1}Q)2^{2j}} \left(1+\|\gamma\|_\infty \right).
\end{align*}
We are going to check that (\ref{hypap}) is satisfied. So we must bound the quantity
$$ \lambda_{k,Q}:= \sum_{j\geq 2} \frac{\mu(2^{j+k+1}Q)}{\mu(2^{k+1}Q)} \alpha_{j,k}  $$
by a constant (independent on $k$ and $Q$).
The coefficient $\alpha_{j,k}$ is estimated by three terms. By using the doubling property for $\mu$, with $N$ large enough we can sum the first term $2^{-N(k+j)}$. For the second term with $N\geq 2$, we use (\ref{gammaint}) to have
\begin{align*}
\sum_{j\geq 2} \frac{\mu(2^{j+k+1}Q)}{\mu(2^{k+1}Q)} \frac{\nu(2^{k}B)}{\nu(2^{k+j+1}B)}  &  2^{-j} \left(\int_{1}^\infty \left(1+2^{j}v^{-2} \right)^{-2N}\gamma(v)^2  v dv \right)^{1/2}  \\
 & \lesssim  \sum_{j\geq 2}  2^{j} \left(\int_{1}^\infty \left(1+2^{j}v^{-2} \right)^{-2N}\gamma(v)^2 v dv \right)^{1/2} \\
 &  \lesssim  \left(\int_{1}^{\infty} v^{2N+1} \gamma(v)^2 dv \right)^{1/2} <\infty.
\end{align*}
For the third term of $\alpha_{j,k}$, we have
 \begin{align*}
\sum_{j\geq 2} \frac{\mu(2^{j+k+1}Q)}{\mu(2^{k+1}Q)} \frac{\mu(2^{k}Q)}{\mu( 2^{j+k}Q)2^{2j}} \left(1+\|\gamma\|_\infty \right) \lesssim \sum_{j\geq 2} 2^{-2j} <\infty.
\end{align*}
We have the desired property due to the additionnal factor $2^{-2j}$, which is obtained by the time-regularity of the semigroup in the case $2-)$. So the assumption $(\ref{hypap})$ is satisfied. \findem

\mb We have the same result for the adjoint operator $T^*$~:
\begin{thm} \label{theo2ap2}  Let $L$ be a generator of a bounded analytic semigroup on $L^2(Y)$. Assume that $(e^{tL^*})_{t>0},(tL^*e^{tL^*})_{t>0}$ and $(t^2L^{2*}e^{tL^*})_{t>0}$ belong to \linebreak[4] ${\mathcal O}_{low}(L^2-L^2) \cap {\mathcal O}_{high}(L^2-L^2)$. Then there exists coefficients $\alpha_{j,k}$ satisfying
\be{hypap2} \Lambda:= \sup_{k\geq 0}
\left[\sum_{j\geq 2} \frac{\mu(2^{j+k+1}Q)}{\mu(2^{k+1}Q)} \alpha_{j,k} \right] <\infty, \ee
such that for all balls $Q\subset X$, for all $k\geq 0,j\geq 2$, for all functions $f$ supported in $S_k(Q)$
$$\left(\frac{1}{\mu(2^{j+k+1}Q)} \int_{S_j(2^kQ)} \left| T^{*}(B_Q^{*}(f))\right|^2 d\mu \right)^{1/2} \leq \alpha_{j,k}(Q) \left(\frac{1}{\mu(2^{k+1}Q)}\int_{S_k(Q)} |f|^2 d\mu \right)^{1/2}.$$
These estimates show with Theorem \ref{theo2h} that $T^{*}$ is $H^1_{\epsilon,mol}-L^1$ bounded for every $\epsilon>0$, with the Hardy space $H^1_{\epsilon, mol}:=H^1_{\epsilon,mol,{({B_Q}^*)}_{Q\in {\mathcal Q}}}$ (which is the Hardy space constructed with the dual operators $B_Q^*$).
\end{thm}

\deme The adjoint operator $T^*$ is given by~:
$$ T^* f(t,x) = \int _{s=t}^{Z}  \left[L^{*} \left(e^{(s-t) L}\right)^{*} f(s,.)\right](x) ds.$$
The parameter $Z$ depends on the time interval $J$, it is defined by~:
$$ Z := \left\{ \begin{array} {ll}
    \infty & \textrm{ if $J=(0,\infty)$} \\
    l & \textrm{ if $J=(0,l)$}
                \end{array} \right. .$$
The argument of the previous theorem can be repeated and we omit details.
\findem

\mb So now we can apply our general result to obtain the following result~:

\begin{thm} Let $L$ be a generator of a bounded analytic semigroup on $L^2(Y)$ such that $(e^{tL})_{t>0},(tLe^{tL})_{t>0}$ and $(t^2L^{2}e^{tL})_{t>0}$ belong to ${\mathcal O}_{low}(L^2-L^2) \cap {\mathcal O}_{high}(L^2-L^2)$ and $(e^{tL^*})_{t>0}$ belongs to ${\mathcal O}_{weak}(L^2-L^{p_0'})$ for a $p_0\in]1,2[$. Then for all $p\in ]p_0,2]$ the operator $T$ is $L^p(X)$-bounded and so we have the maximal regularity on $L^p(Y)$. \\
In addition if we know that $(e^{tL^*})_{t>0},(tL^*e^{tL^*})_{t>0}$ and $(t^2L^{2*}e^{tL^*})_{t>0}$ belong to ${\mathcal O}_{low}(L^2-L^2) \cap {\mathcal O}_{high}(L^2-L^2)$ and $(e^{tL})_{t>0}$ belongs to ${\mathcal O}_{weak}(L^2-L^{q_0})$ for a $q_0\in]2,\infty]$ then for all $p\in [2,q_0[$ the operator $T$ is $L^p(X)$-bounded and so we have the maximal regularity on $L^p(Y)$.
\end{thm}

\deme We use Theorem \ref{theogeneh} for the operator $T$ and $T^*$
whose hypotheses are satisfied thanks to Theorem \ref{theo2ap},
\ref{theo2ap2} and \ref{assumcas}. Then we conclude by duality for
$p\geq 2$. \findem

\begin{rem}  We recall the result of S. Blunck and P.C. Kunstmann (Theorem 1.1 of \cite{BK2}). Assume that the semigroup $e^{tL}$ satisfies for $p_0<2<q_0$ the following $L^{p_0}-L^{q_0}$ off-diagonal estimates for all balls $B\subset Y$~:
\be{BKhy} \left(\frac{1}{\nu(B)} \int_B \left|e^{r_Q^2L}(f)\right|^{q_0} d\nu \right)^{1/q_0} \leq \sum_{j\geq 0} \beta_j \left( \frac{1}{\nu(2^{j+1} B)} \int_{S_j(B)} |f|^{p_0} d\nu \right)^{1/p_0} \ee
with coefficients $\beta_j$ satisfying
\be{BKhyp} \sup_{B\textrm{ ball}} \sum_{j\geq 0} \beta_j \left(\frac{\nu(2^{j+1}B)}{\nu(B)}\right)^{1/p_0} 2^{j(\delta-1)} <\infty. \ee
(We recall that $\delta$ is the homogeneous dimension of the space $Y$.) Then $L$ has maximal regularity on $L^{r}$ for all $r\in(p_0,q_0)$. \\
With our previous theorem, to get maximal regularity for all $r\in (p_0,q_0)$, we need to ask some $L^2-L^2$ off-diagonal estimates and the two following conditions
$$(e^{tL^*})_{t>0} \in {\mathcal O}_{weak}(L^2-L^{p_0'}) \quad \textrm{and} \quad   (e^{tL})_{t>0} \in {\mathcal O}_{weak}(L^2-L^{q_0}). $$
The second condition $(e^{tL})_{t>0} \in {\mathcal O}_{weak}(L^2-L^{q_0})$ implies that we have (\ref{BKhy}) with $\beta_j$ satisfying
$$ \sum_{j\geq 0} \beta_{j} 2^{j} <\infty,$$
which is a weaker condition than (\ref{BKhyp}) for $\delta\geq 2$. So our use of Hardy spaces, adapted to the maximal regularity operator, permits us to show some maximal $L^{r}$ regularity assuming some "stronger" $L^2-L^2$ off-diagonal estimates for the operators $e^{tL}, Le^{tL}, L^2 e^{tL}$ and theirs duals and some "weaker" off-diagonal estimates $L^2-L^{q_0}$ and $L^{p_0}-L^2$ for the operators $e^{tL}$ and $e^{tL^*}$ than those required by the result of S. Blunck and P.C. Kunstmann.
\end{rem}

\mb To finish we will show some results on our Hardy space. First we have the off-diagonal decay (\ref{decay}) and so we are in the particular case of the section \ref{particular}.

\begin{prop} Assume that $(e^{tL})_{t>0} \in {\mathcal O}_{low}(L^2-L^2)$. For $B_Q$ defined by (\ref{bq}) and (\ref{opA}), we have that for all balls $Q \subset X$
\be{decay2} \forall i\geq 0,\ \forall k\geq 0,\ \forall f\in L^2(2^{k}Q), \qquad \left\| B_Q(f)
\right\|_{2,S_i(2^{k}Q)} \leq C 2^{-M''i} \|f\|_{2,2^{k}Q} \ee with an exponent $M''$ as large as we want.
\end{prop}

\deme By definition we have just to prove the decay for the $A_Q$ operator. Let $r$ be the radius of $Q$. As previously, we write $Q=I \times B$ where $I$ is an interval of length $r^2$ and $B$ is a ball in $Y$ of radius $r$. Recall that
$$  A_Q(f)(t,x):=\int_{\sigma=0}^{+\infty} \left[\varphi_{r^2} (t-\sigma)  e^{r^2 L} f(\sigma,.) \right](x) d\sigma.$$
For $i\leq 1$, we just use the $L^2(Y)$-boundedness of $A_Q$ to prove (\ref{decay2}). Then for $i\geq 2$ and $(\sigma,y)\in 2^kQ$ if $(t,x)\in
S_i(2^kQ)$ we have that $d((x,t),(\sigma,y))\simeq 2^{k+i}r$ and by using the definition of the parabolic quasi-distance, we
conclude that either $x\in S_i(2^kB)$ either $t\in S_{2i} (2^{2k}I)$. We will study the two
cases~: \\
First for $x\in S_i(2^k B)$, by the off-diagonal estimate
(\ref{offdiagonallow}) we have the estimate : for all $\sigma>0$
\begin{align*}
 \left\| e^{r^2L}(f(\sigma,.)) \right\|_{2,S_i(2^k B)} \lesssim \frac{\nu(2^kB)}{\nu(2^{i+k}B)} \gamma \left(2^{i+k}\right) \left(\frac{\nu(2^{i+k}B)}{\nu(2^kB)} \right)^{1/2} \|f(\sigma,.)\|_{2,2^{k}B}.
\end{align*}
So by the Minkowski inequality, we obtain
\begin{align*}
\lefteqn{\left\| A_Q(f)(t,.) \right\|_{2,S_i(2^kB)} } & & \\
 & & \lesssim \int_{\sigma=0}^{+\infty} \left(1+ \frac{|t-\sigma|}{r^2} \right)^{-N} \frac{\nu(2^kB)}{\nu(2^{i+k}B)} \gamma\left(2^{i+k}\right) \left(\frac{\nu(2^{i+k}B)}{\nu(2^kB)} \right)^{1/2} \|f(\sigma,.)\|_{2,2^{k}B} \frac{d\sigma}{r^2} \\
  & & \lesssim  \left(\frac{\nu(2^kB)}{\nu(2^{i+k}B)}\right)^{1/2} \gamma\left(2^{i+k}\right) \|f\|_{2,2^{k}Q} \frac{1}{r}.
\end{align*}
Then we integrate for $t\in 2^{2(i+k)}I$ to have
$$ \left\| A_Q(f)\right\|_{2,2^{2(i+k)}I \times S_i(2^kB)} \lesssim \left(\frac{\nu(2^kB)}{\nu(2^{i+k}B)}\right)^{1/2} 2^{i+k} \gamma\left(2^{i+k}\right) \|f\|_{2,2^{k}Q}.$$
For the second case, we have $|t-\sigma|\simeq 2^{2(i+k)}r^2$. By using the $L^2(Y)$-boundedness of the semigroup
$$ \left\| e^{r^2L}(f(\sigma,.)) \right\|_{2,2^{i+k}B} \lesssim  \|f(\sigma,.)\|_{2,2^{k}B}.$$
So by the Minkowski inequality, we obtain
\begin{align*}
  \left\| A_Q(f)(t,.) \right\|_{2,2^{k+i}B)} & \lesssim \int_{\sigma \in 2^{k} I} \left(1+ 2^{2(k+i)} \right)^{-N}  \|f(\sigma,.)\|_{2,2^{k}B} \frac{d\sigma}{r^2} \\
   & \lesssim  2^{-2(k+i)(N-1)} \|f\|_{2,2^{k}Q} \frac{1}{r}.
\end{align*}
So we can conclude that
$$ \left\| A_Q(f)\right\|_{2,S_{2i}(2^{2k}I) \times 2^{i+k}B} \lesssim 2^{-(N-2)(k+i)} \|f\|_{2,2^{k}Q}.$$
With these two cases, we can conclude (for $N$ any large enough integer)
$$ \left\| A_Q(f)\right\|_{2,S_i(2^k Q)} \lesssim \left(2^{-(N-2)i} + \left(\frac{\nu(2^kB)}{\nu(2^{i+k}B)}\right)^{1/2} 2^{i+k} \gamma\left(2^{i+k}\right) \right) \|f\|_{2,2^{k}Q}$$
which with the decay of $\gamma$ permits to prove the result. \findem

\mb With these decay, we have shown that the Hardy spaces
$H^1_{ato}(X)$ and $H^1_{\epsilon, mol}(X)$ are included into the space
$L^1(X)$. In fact we can improve this result, by comparing it with
the classical Hardy space of Coifman-Weiss on $X$.

\begin{prop} \label{aa} Let $\epsilon>0$. The inclusion $H^1_{ato}(X) \subset H^1_{\epsilon, mol}(X) \subset H^1_{CW}(X)$ is equivalent to the fact for all $r>0$, $(e^{rA})^* ({\bf 1}_Y)={\bf 1}_{Y}$ (in the sense of Proposition \ref{inclus}).
\end{prop}

\deme We use the notations of Proposition \ref{inclus}. By using this Proposition, we know that $H^1_{\epsilon, mol}(X) \subset H^1_{CW}(X)$ is equivalent to the fact that for all balls $Q$ of $X$, $A^{*}({\bf 1}_{X}) ={\bf 1}_{X}$ in the sense of $(Mol_{\epsilon,Q})^{*}$. Let $Q=B((t_Q,c_Q),r_Q)$ be fixed. By (\ref{Aetoile}) we know that
$$ A_Q^*(g)(\sigma,x):= \int_{t\in \R^+}  \varphi_{r_Q^2}(t-\sigma) \left[\left(e^{r_Q^2L}\right)^*g(t,.)\right](x) dt.$$
As $\int_\R  \varphi(t) dt =1$, we formally obtain
$$ A_Q^{*}({\bf 1}_X)(\sigma,x) = (e^{r_Q^2L})^{*}({\bf 1}_Y)(x).$$
This equality can be rigorously verified by defining $(e^{r_Q^2L})^{*}({\bf 1}_Y)(x)$ as the continuous linear form on the space
$$ Mol_{\epsilon,r_Q}(Y):= \left\{ f\in L^1(Y),\ \|f\|_{Mol_{\epsilon,r_Q}(Y)}<\infty \right\},$$
where
$$ \|f\|_{Mol_{\epsilon,r_Q}(Y)}:=\sup_{i\geq 0} \|f\|_{2,S_i(Q_Y)} \left(\nu(2^{i}Q_Y)\right)^{1/2} 2^{\epsilon i}.$$
Here we write
$$Q_Y= B(c_Q,r_Q) = \left\{ y\in Y, d_Y(x,c_Q)\leq r_Q \right\}$$
the ball in $Y$. Then the equivalence is a consequence of Proposition \ref{inclus}. \findem

\gb In the paper \cite{ABZ} (which corresponds to the appendix : the chapter 6), the authors have shown that with $-L$ equals to the laplacian on $X$ a complete Riemannian manifold with doubling and Poincaré inequality, the operator $T$ is bounded on $H^1_{CW}(X)$ (not just bounded into $L^1(X)$). This is a better result than the one here because Proposition \ref{aa} applies (see \cite{ABZ}) so
$$ H^1_{ato}(X) \subset H^1_{\epsilon, mol}(X) \subset H^1_{CW}(X) \subset L^1(X). $$
But the $H^1_{CW}$-boundedness is using stronger hypotheses than ours in a specific situation.

\section{Study of the duals of Hardy spaces.} \label{secduality}

In this section, we come back to $(X,d,\mu)$ an abstract space of homogeneous type and we study the dual space of our Hardy spaces. In the Euclidean case the dual space of $H^1_{CW}$ is $BMO$. The other Hardy spaces, discussed in Section \ref{sectionexample}, have duals charaterized by a space of BMO-type (\cite{DY,HM}). In our situation, it is not so easy, and we obtain partial answers. So first, we give the definition of the $Bmo_{\infty}$ space~:

\begin{defn} A function $f\in L^2$ belongs to the space $Bmo_{\infty}$ if
$$ \left\|f \right\|_{BMO} := \sup_{Q \textrm{ ball}} \left(\frac{1}{\mu(Q)}\int_{Q} \left| B_Q^*(f) \right|^2 d\mu \right)^{1/2 } = \left\| M^{\sharp}_2(f) \right\|_{\infty} <\infty.$$
We define $BMO_\infty$ as the completion of $Bmo_\infty$ with this pseudo-norm.
\end{defn}

\mb We note $BMO_\infty$ because the norm is a norm of type $BMO$
and we put the index $\infty$ to note that this space is the closure
of $L^2$. So it could be thought as the space of $BMO$-functions
with some decay at "infinity".

\mb We have the following inclusion~:
\begin{prop} The space $(H^{1}_{ato})^* \cap L^2$ is included into $Bmo_\infty$ and
$$ \forall f\in (H^{1}_{ato})^* \cap L^2 , \qquad \|f\|_{BMO} \leq \|f\|_{(H^{1}_{ato})^*}.$$
We have the same result for the molecular Hardy space $(H^{1}_{\epsilon,mol})^*$.
 \end{prop}

\mb The proof is the same as in Lemma \ref{lem1}.

\mb To characterize the dual space, we should show the other inclusion and we would like to forget the "$L^2$-condition". We give two results in this direction.

\begin{prop} The space $Bmo_\infty$ is included into $(H^{1}_{ato})^* \cap L^2$ and~:
$$ \forall \phi\in Bmo_\infty, \qquad \|\phi\|_{(H^{1}_{ato})^*} \lesssim \|\phi\|_{BMO}.$$
By density, we get that $BMO_\infty \hookrightarrow (H^{1}_{ato})^*$.
\end{prop}

\deme Let $\phi\in Bmo_\infty$ be a function, then $\phi\in L^2$. Let $m=B_Q(f_Q)\in L^2$ be an atom, then we must estimate the quantity~:
$$ \langle \phi,m \rangle = \int \phi(z) B_Q(f_Q)(z) d\mu(z) = \int_{Q} B_Q^*(\phi)(z) f_Q(z) d\mu(z).$$
Then we get~:
$$ \left|\langle \phi,m \rangle \right| \leq \left\|B_Q^*(\phi)\right\|_{2,Q} \left\|f_Q \right\|_{Q} \leq \left\|B_Q^*(\phi)\right\|_{2,Q} \mu(Q)^{-1/2}.$$
Using the $BMO$-norm, we have
$$ \left|\langle \phi,m \rangle \right| \leq  \|\phi\|_{BMO}.$$
So we have shown that $\phi$ can be extended in a continuous linear form on all atoms. By using Lemma \ref{lemdense2}, this proves that $\phi$ is a continuous linear form on $H_{ato}^1$.  \findem

\mb So for the atomic case, we need no other condition to have the two inclusions. For the dual space of molecular space, it seems necessary to have other conditions. Here is an example.

\begin{prop} For a ball $Q$, we write $r_Q$ for the radius. Assume that the operators $A_Q$ depend only on the radius $r_Q$ that is $A_Q=A_{r_Q}$ and $B_Q=B_{r_Q}$.
In this case for all $\epsilon>0$, $Bmo_\infty$ is included into $(H^{1}_{\epsilon,mol})^* \cap L^2$ and~:
$$ \forall \phi\in Bmo_\infty, \qquad \|\phi\|_{(H^{1}_{\epsilon,mol})^*} \lesssim \|\phi\|_{BMO}.$$
Therefore by density, we have $BMO_\infty \hookrightarrow (H^{1}_{\epsilon,mol})^*$.
\end{prop}

\deme Let $m=B_Q(f_Q)$ be an $\epsilon$-molecule. Breaking the integral with the corona
$S_i(Q)$, we have~:
\begin{align*}
 \left|\langle \phi,m \rangle \right| & \leq \sum_{i\geq 0} \int_{S_i(Q)} B_Q^*(\phi)(z) f_Q(z) d\mu(z) \leq \sum_{i\geq 0} \left\| B_Q^*(\phi) \right\|_{2,S_i(Q)}  \left\|f_Q\right\|_{2,S_i(Q)} \\
  & \leq \sum_{i\geq 0} \left\| B_Q^*(\phi) \right\|_{2,S_i(Q)}  \mu(2^{i+1}Q)^{-1/2}2^{-\epsilon i}.
\end{align*}
We want to cover $2^{i+1}Q$ by a union of balls $(\tilde{Q}_k)_{1\leq k \leq \kappa}$ with
$$  \forall k, \quad C^{-1}r_{Q} \leq r_{\tilde{Q}_k} \leq r_Q \textrm{  and  }  \sum_{k=1}^{\kappa} {\bf 1}_{2^{-M}\tilde{Q}_k} \leq 1 $$
where $M$ and $C$ are integers which are large enough and depend only on the space $X$. Let us explain how we can do this. Choose $\left(B(x_i,\frac{1}{3}r_Q)\right)_i$ a maximal collection of disjoint balls. Then we put in the collection $(\tilde{Q}_{k})_{k}$ all the balls $B(x_i,\frac{1}{10} r_Q)$ such that
$$ B(x_i,r_Q) \cap 2^{i+1}Q \neq \emptyset.$$
By maximality, it is easy to see that the collection $(\tilde{Q}_k)$ covers the whole set $2^{i+1} Q$. By disjointness, we have that
$$ \sum_{k=1}^{\kappa}  {\bf 1}_{\frac{1}{3}\tilde{Q}_k} \leq 1.$$
So we can deduce that~:
\begin{align*}
\left\| B_Q^*(\phi) \right\|_{2,S_i(Q)}^2  = \left\| B_{r_Q}^*(\phi) \right\|_{2,S_i(Q)}^2  \lesssim \sum_{k=1}^\kappa \left\| B_{r_Q}^*(\phi) \right\|_{2,\tilde{Q}_k}^2 \lesssim \sum_{k=1}^\kappa \left\| B_{r_Q}^*(\phi) \right\|_{2,\frac{r_Q}{r_{\tilde{Q}_k}}\tilde{Q}_k}^2.
 \end{align*}
Here we use : $r_{Q}\geq r_{\tilde{Q}_{k}}$. In addition the radius of the ball $\frac{r_Q}{r_{\tilde{Q}_k}}\tilde{Q}_k$ is now equal to the radius $r_Q$, also
 $$ \left\| B_{r_Q}^*(\phi) \right\|_{2,\frac{r_Q}{r_{Q_k}}\tilde{Q}_k}^2 = \left\| B_{\frac{r_Q}{r_{\tilde{Q}_k}}\tilde{Q}_k}^*(\phi) \right\|_{2,\frac{r_Q}{r_{\tilde{Q}_k}}\tilde{Q}_k}^2 \lesssim \mu(\frac{r_Q}{r_{\tilde{Q}_k}}\tilde{Q}_k) \|\phi\|_{BMO}.$$
By using the doubling property of the measure and the fact that $r_Q\simeq r_{\tilde{Q}_k}$, we get
\begin{align*}
\left\| B_Q^*(\phi) \right\|_{2,S_i(Q)}^2 & \lesssim \sum_{k=1}^\kappa \|\phi\|_{BMO}^2 \mu(\tilde{Q}_k) \lesssim \|\phi\|_{BMO}^2 \sum_{k=1}^\kappa \mu(3^{-1}\tilde{Q}_k) \\
 & \lesssim \|\phi\|_{BMO}^2  \mu(2^{i+1}Q).
\end{align*}
Therefore we have the estimate
$$ \left|\langle \phi,m \rangle \right| \lesssim \sum_{i\geq 0}  \|\phi\|_{BMO} 2^{-\epsilon i} \lesssim \|\phi\|_{BMO}.$$
Thus we have shown that $\phi$ is linearly continuous on all the
$\epsilon$-molecule. With the help of Lemma \ref{lemdense2}, we can
extend it to a continously linear form on the whole molecular space.
\findem

\begin{rem} Assume that $A_Q$ is defined by the help of a semi-group $A_Q=A_{r_Q}=e^{-r_QL}$ for a generator $L$.
Then the previous condition is satisfied and with some good
conditions about $L$, we can characterize the dual space
$(H^1_{\epsilon,mol})^* \cap L^2$ as $Bmo_\infty$. It is interesting
to notice that this space does not depend on $\epsilon$.
\end{rem}

\gb Here we have not found a general answer for the dual space
$(H^1_{ato})^*$. We have only the result $Bmo_\infty=L^2 \cap
(H^1_{ato})^*$. Now to have a complete result for the duality, we
must forget the ``$L^2$ condition''. In such an abstract case, it
seems very difficult to do this. The difficulty is to have a
representation of a linearly continuous form of $H^1_{ato}$. Let $l$
belongs to $(H^1_{ato})^*$, by the definition of atom we have that
the operator $B_Q$ is continuous from $L^2(Q)$ into $H^1_{ato}$. So
we can compute $l \circ B_Q$, which is continuous from $L^2(Q)$ into
$\R$. By the Riesz representation theorem, we know that there exists
$h_Q\in L^2(Q)$ such that for all $f\in L^2(Q)$
$$  l\circ B_Q( f) := \int_Q h_Q(x) f(x) d\mu(x).$$
Now to have a good representation of $l$, we ``need to invert'' the operator $B_Q$ to have
formally~: for all atom $g$ associated to the ball $Q$
$$  l(g) = \int_Q h_Q(x) B_Q^{-1}g(x) d\mu(x) = \int_Q B_Q^{-1*}h_Q(x) g(x) d\mu(x). $$
Then we need to define a ``good'' function $\phi$ (which does not depend on $Q$) such that for each ball $Q$ and for each atom $g$ associated to $Q$ we have
 $$  l(g) = \int_Q B_Q^{-1*}h_Q(x) g(x) d\mu(x) = \int_X \phi(x) g(x) d\mu(x). $$
Here if we do not want to invert, we need a decomposition of identity with the $B_Q$ operators. That is why, in \cite{DY,HM} the authors use a Calder\'on reproducing formula to identify their dual space. \\
So there are two questions to solve and we do not know how to do this ; other informations on the collection of operators ${\mathbb B}=(B_Q)_{Q\in {\mathcal Q}}$ seem to be necessary, but we do not know at this time which ones.

\gb We have seen that the space $(H^1_{ato})^{*}$ is probably too big to be identified with a
 BMO-space. However we are going to show that the subspace $L^2 \cap (H^1_{ato})^{*}$ is dense in the whole space $(H^1_{ato})^{*}$ for a weak topology. Let us study the topology of $Bmo_\infty$ in $(H^1_{ato})^{*}$ for the weak *
topology. We recall that $S$ is the set of functions $f\in
H^1_{ato}$ so that there exists a finite decomposition into atoms
$(m_i)$ satisfying~:
$$ f = \sum_{i=1}^{n} \lambda_i m_i \qquad \textrm{with} \qquad \|f\|_{H^1_{ato}} \geq 10^{-1} \left(\sum_i |\lambda_i|\right).$$
By Lemma \ref{lemdense2}, we know that $S$ is dense in $H^1_{ato}$.

\begin{prop} \label{thdualite} If the space $Bmo_\infty$ is dense in $L^2$ (for the strong topology) then $Bmo_\infty$ is dense in $(H^1_{ato})^*$ for the weak * topology of $S^*$.
\end{prop}

\deme  We claim that the space $Bmo_\infty$ is total, that is~:
\be{total} \left\{ f\in S;\ \forall \phi\in Bmo_\infty,\ \langle f,
\phi \rangle =0 \right\}=\left\{0\right\}. \ee To prove this, let
$f\in S$ be a function in the left set. Then $f$ has a finite atomic
decomposition so it is an $L^2$-function, and we have
$$  \forall \phi\in Bmo_\infty\subset L^2, \qquad _{L^2} \langle f , \phi \rangle_{L^2} =0.$$
As $Bmo_\infty$ is assumed to be dense in $L^2$ and $f$ belongs to $L^2$, we can deduce that $f=0$, which proves (\ref{total}). \\
We use a general fact : Theorem 4 of \cite{dix} to obtain the density for the weak * topology. \findem

\mb We can have a more precise theorem in the case of the section \ref{particular}~:

\begin{prop} Assume that the assumptions (\ref{decay}) of the section \ref{particular} are satisfied. The space $L^\infty\cap L^2$ is dense in $(H^1_{ato})^*$ for the weak * topology of $(H^1_{ato})^*$. For $\epsilon>0$, $L^\infty\cap L^2$ is dense in $(H^1_{\epsilon,mol})^*$ for the weak * topology.
\end{prop}

\deme The proof is the same as the one of Proposition \ref{thdualite}. We will deal with the case $\epsilon\in]0,\infty]$ and we will prove the two claims.
The fact that the space $L^\infty$ is total, means that~:
\be{total2}
\left\{ f\in H^1_{\epsilon, mol},\ \forall \phi\in L^\infty\cap L^2 \subset (H^1_{\epsilon, mol})^{*},\ \langle f, \phi \rangle =0 \right\}=\left\{0\right\}. \ee
This fact is obvious because the function $f\in L^1$ (due to Proposition \ref{contL1}).
As for Proposition \ref{thdualite}, we use the Theorem 4 of  \cite{dix} to conclude. \findem

\mb In fact to have a complete representation theorem for the dual space, we probably need to
 make some new assumptions. In \cite{DY,DY1,HM}, the authors  characterize the
 dual space by a BMO-space, by using an equivalent definition of their Hardy spaces with
 tent spaces.Using molecular decomposition in tent spaces, they obtain some molecular
 decomposition of their Hardy spaces. Without other assumptions, our molecular decomposition is strictly more restrictive than theirs. So in the general case, we think that the dual space of our Hardy spaces is bigger than a BMO-space. We have seen that our Hardy spaces are "big" enough to obtain a good interpolation result with the scale of Lebesgue spaces, but they seem to be too "small" to have a fine dual space.

\vspace{0.1in} {\bf Acknowledgements :} The authors are indebted to
professor Pascal Auscher for suggesting the topic, numerous
discussions and advices to improve this paper. The second author
would like to express many thanks to Department of Mathematics,
Paris-Sud University for its hospitality.

\end{document}